\documentclass{amsart}
\usepackage[italian,english]{babel}
\usepackage{amsaddr}

\usepackage[utf8]{inputenc}

\usepackage{enumitem}

\usepackage{biblatex}
\usepackage{csquotes}
\usepackage{lmodern}

\usepackage{amsmath}

\usepackage{amssymb}

\usepackage{bbm}

\usepackage{fancyhdr}

\numberwithin{equation}{section}

\usepackage{xcolor}
\usepackage{geometry}
\newtheoremstyle{theorem2}
{8pt}
{8pt}
{\itshape}
{}
{\bfseries}
{.}
{.5em}
{}

\newtheoremstyle{definition2}
{8pt}
{8pt}
{\itshape}
{}
{\bfseries}
{.}
{.5em}
{}

\theoremstyle{theorem2}

\newtheorem{lemma}{Lemma}[section]
\newtheorem{teo}[lemma]{Theorem}

\newtheorem{cor}[lemma]{Corollary}
\newtheorem{prop}[lemma]{Proposition}

\theoremstyle{definition2}

\newtheorem{deff}[lemma]{Definition}
\newtheorem{Remark}[lemma]{Remark}

\begin{document}

\title{Transport densities and congested optimal transport problem in the Heisenberg group}

\author{Michele Circelli, Giovanna Citti} 
\address{University of Bologna, Department of Mathematics, 40126 Bologna, Italy}

\keywords{Congested optimal transport, Wardrop equilibrium, Heisenberg group, transport density} 

\begin{abstract}
We adapt the problem of continuous congested optimal transport to the Heisenberg group, equipped with a sub-Riemannian metric. Originally introduced in the Euclidean setting by Carlier, Jimenez, and Santambrogio as a path-dependent variant of the Monge-Kantorovich problem, we significantly restrict the set of admissible curves to horizontal ones. We establish the existence of equilibrium configurations as solutions to a convex minimization problem over a suitable set of measures on horizontal curves. This result is achieved through the notions of horizontal transport density and horizontal traffic intensity.
\end{abstract}

\maketitle

\section{Introduction}
The first formalization of the concept of congested optimal transport in a discrete, Euclidean setting can be traced back to 1952, when Wardrop introduced a notion of equilibrium governing traffic congestion on a network in \cite{Wardrop}. In 1956, Beckmann et al. proposed a variational characterization for such equilibria in \cite{Beckmann1}. For a comprehensive overview, see \cite{Santambrogiolibro} or \cite{Santambrogio2}. In a slightly different context, Bouchitté and Buttazzo introduced the notion of transport density in \cite{Buttazzo1} and \cite{Buttazzo2}, which measures the amount of transport occurring along geodesics within a given region for an optimal transport plan. In \cite{Morel}, Bernot et al. introduced a similar but more abstract framework for describing traffic congestion in networks, expressed in terms of traffic plans, which are probability measures over all possible curves between origin-destination pairs. A few years later, Carlier et al. (in \cite{Santambrogio1}) used this optimal transport approach to propose a continuous version of the traffic congestion problem in networks. They proved the existence of Wardrop equilibria through a convex minimization problem over the set of traffic plans, which is the continuous counterpart of the one proposed in \cite{Beckmann1} for the discrete setting. The traffic intensity they introduced can be considered as the reformulation, in the abstract framework introduced in \cite{Morel}, of notion of transport density. It is worth noting that this minimization problem allows for duality and numerical simulations, as discussed in \cite{Peyre2} and \cite{Peyre}. However, this minimization problem over a set of measures over curves presents many technical difficulties. To address this, in \cite{Brasco} Brasco et al. proved  the equivalence with a technically simpler minimization problem over a set of vector fields with prescribed divergence, introduced in \cite{Beckmann2}. They also proved that Wardrop equilibria are supported on the integral curves, in the DiPerna-Lions sense, of optimal vector fields. For a regularity result about optimal vector fields, see \cite{Colombo}. This new formulation is easier to handle and admits a dual formulation as a classical problem in the calculus of variations. Optimality conditions for this formulation can be written as a PDE, which is the $q$-Laplacian in the simplest case but becomes quite degenerate in more realistic congestion models. See also \cite{brasco2013congested} and \cite{Santambrogioregularity}.
\\

In this paper, we address the continuous congested optimal transport problem in the Heisenberg group $\mathbb{H}^n$, the simplest non-commutative Lie group, which naturally arises in the description of the phase space. This setting is characterized by the choice of $2n$ vector fields and a metric on the sub-bundle of the tangent bundle generated by them. Since the number of given vector fields is strictly smaller than the dimension of the space, the metric is non-Riemannian at every point. All intrinsic objects of the space are defined in terms of these vector fields: in particular, displacement occurs only along their integral curves, called horizontal curves, and the presence of this metric naturally leads to the definition of a sub-Riemannian distance. The Heisenberg group was also proposed by Petitot and Tondut in \cite{Petitot} to describe the functional geometry of the visual cortex (see also \cite{cittisarti}, where the problem was expressed using the tools of sub-Riemannian geometry): they describe the propagation of signals along horizontal curves. Additionally, in \cite{Misic}, the authors state that the propagation of visual signals gives rise to congestion phenomena. With this in mind, we adapt to the Heisenberg group results contained in \cite{Santambrogio1}, which were specific to the Euclidean setting. Addressing this problem in the Heisenberg group is very interesting from a mathematical viewpoint due to its close link with optimal transport theory in this setting. In this direction, the first significant work is \cite{Rigot}, where the authors demonstrated the existence of solutions to the Monge problem associated with the square of the sub-Riemannian distance and established a Brenier-McCann representation theorem. This result was further generalized in \cite{FigalliRifford} to more general (non-homogeneous) sub-Riemannian structures. The first existence result for solutions to the Monge problem associated with the sub-Riemannian distance was provided by De Pascale and Rigot in \cite{DePascale2}. For similar results in more general metric spaces, see \cite{Cavalletti3} and \cite{cavalletti2018overview}.

In this context, we introduce the notion of horizontal transport density, adapting to the geometry of the Heisenberg group the notion previously introduced in the Euclidean setting in \cite{Buttazzo1} and \cite{Buttazzo2}. Following the approach in \cite{Santambrogio}, we establish some summability results of the horizontal transport density for a range of $p$ depending on the geodetic homogeneity of the space. It is worthwhile to note that in $\mathbb{H}^n$, the geodesic dimension does not coincide with either the topological or the homogeneous dimension (see \cite{Juillet}). These results will be instrumental in defining the congested optimal transport problem in this setting. To address this problem, the first obstacle to overcome is lack of non-trivial geodesic subsets in the Heisenberg group, as discussed in \cite{Monti2}. Then, we define a horizontal traffic plan as a probability measure over the space of horizontal curves. The associated horizontal traffic intensity, which measures the amount of traffic along horizontal curves, generalizes the concept of horizontal transport density. Subsequently, inspired by \cite{Santambrogio1}, we introduce a convex minimization problem over the set of horizontal traffic plans, whose solutions are Wardrop equilibria in $\mathbb{H}^n$. To prove the existence of such solutions, we use results from \cite{Santambrogio1} that utilize abstract measure theory instruments. As optimality conditions, we find that such solutions are equilibrium configurations. In fact, traffic plans that solve the minimization problem are supported on geodesics with respect to a suitable congested metric (depending on the traffic plan itself through the traffic intensity). Moreover, they induce transport plans that solve a Monge-Kantorovich problem associated with the congested metric. One of the delicate aspects of the paper is defining this metric, due to the difficulty in defining a length induced by $p$-summable weights. 

The problem addressed in this paper, in addition to its intrinsic interest - being an optimal transport problem in a sub-Riemannian setting that has not yet been explored in the literature - has another important aspect: it admits two additional formulations. The first is a minimal-flow-type formulation, expressed as a minimization problem over the space of $L^p$-summable horizontal vector fields. The second is a dual formulation, which takes the form of a classical problem in the calculus of variations. For both these formulations see \cite{circelli2024continuous}. The latter is particularly relevant for the study of regularity properties of solutions to the associated partial differential equations in the Heisenberg group. For instance, the Lipschitz regularity for solutions to an orthotropic equation arising from this model is established in \cite{circelli2025lipschitz}, while \cite{circelli2025gradientestimates} addresses the corresponding parabolic case.
\\

The paper is organized as follows: in Section 2, we introduce the Heisenberg Group and recall some well-known results on optimal transport theory in this setting. In Section 3, we introduce the notion of transport density and provide conditions that ensure its $L^p$-summability. In Section 4, we demonstrate that an $L^p$ function induces a weighted distance. In Section 5, we adapt the problem of congested optimal transport to the Heisenberg group and prove the existence of Wardrop equilibria in this setting.

\section{Preliminaries}
\subsection{The Heisenberg group $\mathbb{H}^n$}
Let $n\geq1$. The $n$-th \textit{Heisenberg group} $\mathbb{H}^n$ is the connected and simply connected Lie group, whose Lie algebra $\mathfrak{h}^n$ is stratified of step $2$, i.e. 
$$ \mathfrak{h}^n= \mathfrak{h}^n_1 \oplus \mathfrak{h}^n_2,$$
where $\mathfrak{h}_1= \mathrm{span} \{X_1, \dots, X_n, X_{n+1}, \dots, X_{2n}\}$ is called \textit{horizontal layer}, $\mathfrak{h}_2=\mathrm{span} \{X_{2n+1} \}$ and the only non-trivial bracket-relation between the vector fields $X_1, \dots, X_{2n+1}$ is 
\begin{equation*}
	[X_j, X_{n+j}]=X_{2n+1},\quad \forall j=1,\ldots n.
\end{equation*}

The horizontal layer $\mathfrak{h}^n_1$ induces a sub-bundle of the tangent bundle, that we denote by $H \mathbb{H}^n$ and whose fibre at any $q \in \mathbb{H}^n$ is 
$$H_q \mathbb{H}^n=\mathrm{span} \{ X_1(q), \dots, X_n(q), X_{n+1}(q), \dots, X_{2n}(q) \}.$$

We call it \textit{horizontal bundle} and any section \textit{horizontal vector field}. 

For any $q \in \mathbb{H}^n$, we consider on $H_q \mathbb{H}^n$ an inner product $\langle \cdot, \cdot \rangle_{H,q}$ that makes $\{X_1(q), \dots, X_{2n}(q)\}$ an orthonormal basis. We denote by $|\cdot|_{H,q}$ the norm associated with such an inner product. 
As it is common in Riemannian geometry, we drop the index $q$ in the inner product, writing $\langle \cdot, \cdot \rangle_H,\forall q\in\mathbb{H}^n$. The same convention shall be adopted for the norm.

Given a smooth vector field $X:\mathbb{H}^n\rightarrow H\mathbb{H}^n$, $t\in(-\varepsilon,\varepsilon)$ and $q_0\in\mathbb{H}^n$, we denote by $\exp(tX)(q_0):=\sigma(t)$, where $\sigma$ is the curve that solves 
\begin{equation*}
	\begin{cases}
		\dot{\sigma}(t)=X(\sigma(t)),\\
		\sigma(0)=q_0.
	\end{cases}
\end{equation*}
If $X\in\mathfrak{h}^n$ and $q_0\in\mathbb{H}^n$, then the previous map is well-defined for any $t\in\mathbb{R}$, see for instance \cite[Proposition 2.1.53]{Bonfiglioli}. Moreover it holds that $$\exp(tX)(q_0)=q_0\cdot\exp(tX)(e),\quad \forall q_0\in\mathbb{H}^n,$$
where $\cdot$ denotes the group law and $e$ is the identity element.

Since the group is also nilpotent, then the exponential map $\mathrm{exp}: \mathfrak{h}^n\to \mathbb{H}^n$ defined as
\begin{equation*}
	\exp(X):=\exp(X)(e)
\end{equation*}
is a global diffeomorphism. Hence, every $q\in\mathbb{H}^n$ can be written in an unique way as
\begin{equation*}
	\label{coordinates}
	q= \mathrm{exp}(x_1 X_1+\dots+x_n X_n+x_{n+1}X_{n+1}+\dots x_{2n}X_{2n}+x_{2n+1}X_{2n+1}),
\end{equation*} 
where $x_i\in\mathbb{R},\forall i=1,\ldots,2n+1$. 
This induces a system of globally defined coordinates on $\mathbb{H}^n$, by identifying any point $q \in \mathbb{H}^n$ with $(x_1,\ldots,x_{2n+1})\in \mathbb{R}^{2n+1}$. From now on we work in this system of coordinates and we write $x=(x_1,\ldots,x_{2n+1})\in\mathbb{H}^n\simeq\mathbb{R}^{2n+1}$. 

In this system of coordinates the group law reads as   
\begin{equation*}\label{grouplaw}
	x \cdot y := \bigg{(}x_1+y_1,\ldots,x_{2n}+y_{2n}, x_{2n+1}+y_{2n+1}+\frac{1}{2}\sum_{j=1}^{2n} (x_jy_{n+j}-x_{n+j}y_j)\bigg{)},
\end{equation*}
and the left-invariant vector fields $X_1,\ldots,X_{2n+1}$ as
\begin{equation*}
	\begin{cases}
		X_j := \partial_{x_j} -\frac{x_{n+j}}{2} \partial_{x_{2n+1}} \,, j=1,\dots,n,\\
		X_{n+j} := \partial_{x_{n+j}}+\frac{x_j}{2} \partial_{x_{2n+1}},\ \ j=1,\dots,n,\\
		X_{2n+1}=\partial_{x_{2n+1}}.	
	\end{cases}
\end{equation*}

Always in coordinates, the unit element $e\in\mathbb{H}^n$ is $0_{\mathbb{R}^{2n+1}}$, the center of the group is 
\begin{equation*}
	L := \{(0,\ldots,0,x_{2n+1}) \in \mathbb{H}^n ;\; x_{2n+1} \in \mathbb{R}\},
\end{equation*}
and, according to the two steps stratification of $\mathfrak{h}^n$, there exists family of non-isotropic dilations that reads as
\begin{equation*}
	\delta_\lambda((x_1,\ldots,x_{2n+1})):= (\lambda x_1,\ldots,\lambda x_{2n},\lambda^2x_{2n+1}),\quad \forall x\in\mathbb{H}^n, \forall \lambda>0.
\end{equation*} 

The homogeneous dimension of $\mathbb{H}^n$ is 
\begin{equation}\label{homogdim}
	N:=\sum_{j=1}^2j\dim(\mathfrak{h}^n_j)=2n+2.
\end{equation}

See \cite{Bonfiglioli}, in particular Section 2, for a general overview on Carnot groups.

Let us just remark that the Lebesgue measure $\mathcal{L}^{2n+1}$, which we shall also denote by $dx$, is the Haar measure of the group $\mathbb{H}^n\simeq\mathbb{R}^{2n+1}$.

Let now consider an open set $\Omega\subseteq\mathbb{H}^n$ and a measurable function $f:\Omega\to\mathbb{R}$. We denote by $$\nabla_Hf=(X_1f,\ldots,X_{2n}f),$$
where $X_jf$ is the derivative of $f$ in the horizontal direction $X_j$, in the sense of distributions. 

Let $1\leq p\leq\infty$, then the space
\begin{equation}\label{horsob}
    HW^{1,p}(\Omega):=\left\{f:\Omega\to\mathbb{R} \text{ measurable : }f\in L^{p}(\Omega),\nabla_Hf\in L^{p}(\Omega)\right\},
\end{equation}
equipped with the norm
\begin{equation*}
    \|f\|_{HW^{1,p}(\Omega)}:=\|f\|_{L^p(\Omega)}+\|\nabla_H f\|_{L^p(\Omega)},
\end{equation*}
is a Banach space. Moreover, for $1\leq p<\infty$, we denote by
\begin{equation*}
    HW_0^{1,p}(\Omega):=\overline{C_0^{\infty}(\Omega)}^{HW^{1,p}(\Omega)}
\end{equation*}
and by
\begin{equation*}
    HW^{-1,p'}(\Omega):=( HW_0^{1,p}(\Omega))'.
\end{equation*}


\subsubsection{\textbf{Carnot-Carath\'eodory distance}}
We can equip $\mathbb{H}^n$ with a sub-Riemannian distance, also known as Carnot-Carath\'eodory distance, that makes it a polish space.

We say that a Lipschitz curve $\sigma\in \text{Lip}([a,b],\mathbb{R}^{2n+1})$, is \textit{horizontal} if its velocity vector $\dot{\sigma}(t)$ belongs to $H_{\sigma(t)}\mathbb{H}^n$ at almost every $t\in[a,b]$ where it exists.
We will denote by
\begin{equation*}
	H([a,b],\mathbb{R}^{2n+1}):=\left\{\sigma\in \text{Lip}\left([a,b],\mathbb{R}^{2n+1}\right):\sigma\text{ is horizontal}\right\}.
\end{equation*}

Given a non negative continuous function $\phi\in C(\mathbb{H}^n,\mathbb{R}_+)$, the \textit{horizontal length of} $\sigma\in H([a,b],\mathbb{R}^{2n+1})$ \textit{weighted by} $\phi$ is
\begin{equation}\label{3agosto}
    L_\phi(\sigma):=\int_a^b\phi(\sigma(t))|\dot\sigma(t)|_Hdt.
\end{equation}
When $\phi\equiv1$, it reads as the \textit{horizontal length} of $\sigma$
\begin{equation}\label{horlength}
\ell_H(\sigma):=\int_a^b|\dot{\sigma}(t)|_H dt.
\end{equation}
Given $x,y\in\mathbb{H}^n$, one can define the \textit{Carnot-Caratheodory distance} (shortly \textit{CC-distance}) between them as 
\begin{equation}\label{CCdistance}
	d_{CC}(x,y):= \inf \  \left\{  \ell_H( \sigma) \ | \  \sigma\in H\left([a,b],\mathbb{R}^{2n+1}\right) , \ \sigma(a)=x, \ \sigma(b)=y \right\}.
\end{equation}

The Chow-Rashevskii theorem guarantees that this distance is well-defined and it induces the Euclidean topology on $\mathbb{H}^n\simeq\mathbb{R}^{2n+1}$, see for instance \cite[Theorem 3.31]{barilariagrachev}. 

This distance is left invariant and 1-homogeneous with respect to the dilations, i.e.
\begin{equation*}
	d_{CC}(x \cdot y, x\cdot z) = d_{CC}(y,z) \quad \text{and} \quad d_{CC}(\delta_\lambda(y), \delta_\lambda(z)) = \lambda\,d_{CC}(y,z)
\end{equation*}
for all $x$, $y$, $z\in\mathbb{H}^n$ and all $\lambda>0$. 

\subsubsection{\textbf{Geodesics}}
We call \textit{minimizing horizontal curve} any $\sigma\in H([a,b],\mathbb{R}^{2n+1})$ such that
\begin{equation*}
	\ell_H(\sigma)=d_{CC}(\sigma(a),\sigma(b)).
\end{equation*}
According to the terminology used in literature, we call \textit{geodesic} any minimizing horizontal curve parametrized proportionally to the arc-length, i.e. $d_{CC} (\sigma(t),\sigma(t'))=|t-t'|v$, where $v=\frac{d_{CC}(\sigma(a),\sigma(b))}{b-a}$ is the constant speed of $\sigma$.  The sub-Riemannian version of the Hopf-Rinow theorem implies that $(\mathbb{H}^n, d_{CC})$ is a geodesic space, see for instance \cite[Theorem 2.4]{FigalliRifford}.

In the Heisenberg group geodesics have been computed explicitly and so it was possible to detect a set in which these are unique. See \cite{Rigot}, \cite{DePascale2} or \cite{Monti}.

We denote by
\begin{equation} \label{KAPPA}
	E:= \{(x,y)\in \mathbb{H}^n\times \mathbb{H}^n;\,\, x^{-1}\cdot y \not \in L\},
\end{equation}
then it holds the following characterization for geodesics parametrized on $[0,1]$.

\begin{teo}\label{geod} 
A non-trivial geodesic starting from $0$ is the restriction to $[0,1]$ of the curve $$\sigma_{\chi,\theta}(t)=\left(x_1(t),\ldots,x_{2n+1}(t)\right)$$ either of the form
\begin{align}\label{geodform}
    &x_j(t)=\frac{\chi_j\sin(\theta t)-\chi_{n+j}\left(1-\cos(\theta t)\right)}{\theta},\quad j=1,\ldots,n\\
    &x_{n+j}(t)=\frac{\chi_{n+j}\sin(\theta t)+\chi_{j}\left(1-\cos(\theta t)\right)}{\theta},\quad j=1,\ldots,n\\
    &x_{2n+1}(t)=\frac{|\chi|^2}{2\theta^2}\left(\theta t-\sin(\theta t)\right),
\end{align}
for some $\chi \in \mathbb{R}^{2n}\setminus\{0\}$ and $\theta\in [-2\pi,2\pi]\setminus\left\{0\right\}$, or of the form
\begin{equation*}
    \left(x_1(t),\ldots,x_{2n+1}(t)\right)=\left(\chi_1t,\ldots,\chi_{2n}t,0\right), 
\end{equation*}
for some $\chi \in \mathbb{R}^{2n}\setminus\{0\}$ and $\theta=0$. In particular, it holds $$|\chi|_{\mathbb{R}^{2n}}=|\dot\sigma|_H=d_{CC}(0,\sigma(1)).$$
	
In particular it is a horizontal curve and it holds:
\begin{enumerate}
    \item For all $(x,y)\in E$, there is a unique geodesic $x\cdot\sigma_{\chi,\theta}$ parametrized on $[0,1]$ between $x$ and $y$, for some $\chi \in \mathbb{R}^{2n}\setminus\{0\}$ and some $\varphi\in(-2\pi,2\pi)$.
    \item If $(x,y)\not\in E$, then $x^{-1}\cdot y = (0,\ldots,0,z_{2n+1})$ for some $z_{2n+1}\in \mathbb{R}\setminus\{0\}$. Hence, there are infinitely many geodesics parametrized on $[0,1]$ between $x$ and $y$: they are all the curves of the form $x\cdot \sigma_{\chi,2\pi}$, if $z_{2n+1}>0$,  $x\cdot \sigma_{\chi,-2\pi}$, if $z_{2n+1}<0$, for any $\chi \in \mathbb{R}^{2n}$ such that $|\chi|_{\mathbb{R}^{2n}} = \sqrt{4\pi|z_{2n+1}|} $.
\end{enumerate}
\end{teo}

From the previous theorem follows that $(\mathbb{H}^n, d_{CC})$ is a non-branching metric space, then any two geodesics that coincide on a non-trivial interval coincide on the whole intersection of their intervals of definition.

We denote by 
\begin{equation*}
    \textnormal{Geo}(\mathbb{H}^n):=\left\{\sigma\in H\left([0,1],\mathbb{R}^{2n+1}\right):\sigma\textnormal{ is a geodesic}\right\}.
\end{equation*}

From the theory of Souslin sets and general theorems about measurable selections, it follows that there exists a  map 
\begin{equation}\label{19marzo1}
	S:\mathbb{H}^n\times\mathbb{H}^n\rightarrow\textnormal{Geo}(\mathbb{H}^n)
\end{equation}
such that for every $x,y \in \mathbb{H}^n$, the value  $S(x,y)=\sigma_{x,y}$ is a geodesic between $x$ and $y$ and $S$ is $\gamma$-measurable for any positive Borel measure $\gamma$ on $\mathbb{H}^n\times\mathbb{H}^n$ (see for instance \cite[Theorem 6.9.2 and Theorem 7.4.1]{Bogachev}). Moreover $S$ is continuous, hence Borel, on the set $E$, defined in \eqref{KAPPA}. The map $S$ is often called \textit{selection of geodesics} map. 

If $e_t$ is the evaluation map at $t\in[0,1]$, i.e. $e_t(\sigma) := \sigma(t)$ for all $\sigma \in C\left([0,1],\mathbb{R}^{2n+1}\right)$, then the map
\begin{equation}\label{geod01}
	S_t:=e_t\circ S:\mathbb{H}^n\times\mathbb{H}^n\rightarrow\mathbb{H}^n,
\end{equation}
associates with any two points $x,y\in\mathbb{H}^n$, the point $S_t(x,y)$ of $\mathbb{H}^n$ at distance $t \, d_{CC}(x,y)$ from $x$, on the selected geodesic $S(x,y)$ between $x$ and $y$. If we fix $\overline{y}\in\mathbb{H}^n$ and $t\in(0,1)$, then the function $S_t(\cdot,\overline{y})$ is $C^{\infty}$ on $\mathbb{H}^n\setminus(\overline{y}\cdot L)$ and it holds that
\begin{equation}\label{det}
	\det D_x(S_t(x,\overline{y}))\geq(1-t)^{2n+3},
\end{equation}
for all $x\in\mathbb{H}^n\setminus(\overline{y}\cdot L)$. Moreover, for any $y\in\mathbb{H}^n$ and any $A\subset \mathbb{H}^n\setminus \left(y\cdot L\right)$,
\begin{equation}\label{MCP}
	\mathcal{L}^{2n+1}(A) \leq \frac{1}{(1-t)^{2n+3}} \,\mathcal{L}^{2n+1}(S_t(A,y)),
\end{equation}
which means that $(\mathbb{H}^n,d_{CC},\mathcal{L}^{2n+1})$ satisfies a so-called Measure Contraction Property $MCP(0,2n+3)$: see \cite[Section 2]{Juillet} for proofs of these results. The measure contraction property is a generalization to metric measure spaces of the concept of Ricci curvature bounded by below. This notion was introduced by Otha in \cite{ohta2007measure}: it controls the distortion of measures along geodesics. Recently in \cite{badreddine2020measure} the authors proved that every two-step compact sub-Riemannian manifold and every Lipschitz Carnot group satisfy $MCP(0,R)$, for some $R>0$. See also \cite{barilari2018sharp}, \cite{rifford2013ricci}, \cite{rizzi2016measure} and references therein for further results in this direction.

We explicitly recall that formulas \eqref{det} and \eqref{MCP} extend to this setting 
relations well known in the Euclidean setting, where quantity $1-t$ is raised to  dimension of the space. Here the exponent $2n+3$ is neither the topological dimension  $2n+1$, nor the homogeneous dimension $N=2n+2$ of the Heisenberg group. It is the so called geodesic dimension of the space and  Juillet, in the \cite[Remark 2.3]{Juillet}, shows that this exponent is sharp.

\subsection{Optimal transport in $\mathbb{H}^n$}\label{MonKant}

This section contains some well-known results about optimal transportation theory in the Heisenberg Group. Precisely we collect some results about the Monge-Kantorovich problem associated with Carnot-Carath\'eodory distance, 
following the presentation in \cite{DePascale2}.

Given a Polish metric space $(M,d)$, we denote by $\mathcal{M}_+(M)$ the set of positive and finite Radon measures on $M$;  we denote by $\mathcal{P}(M)$, resp. $\mathcal{P}_c(M)$, the subset of probability measures on $M$, resp. the subset of probability measures on $M$ with compact support. 

If $(M_1,d_1)$ and $(M_2,d_2)$ are two Polish metric spaces, $f:M_1\to M_2$ is a Borel map and $\mu\in\mathcal{M}_+(M_1)$, we denote by $f_{\#}\mu\in\mathcal{M}_+(M_2)$ the measure defined as $f_{\#}\mu(A):=\mu(f^{-1}(A))$, $\forall A$ Borel set in $M_2$.
Let $\mu,\nu\in\mathcal{P}_c(\mathbb{H}^n)$, we denote by
\begin{equation*}
	\Pi(\mu,\nu)=\big{\{}\gamma\in\mathcal{P}(\mathbb{H}^n\times\mathbb{H}^n): (\pi_1)_{\#}\gamma=\mu, (\pi_2)_{\#}\gamma=\nu\big{\}}
\end{equation*}
the set of transport plans between $\mu$ and $\nu$, where $\pi_1$ and $\pi_2$ are the projection on the first and second factor, respectively. The set $\Pi(\mu,\nu)$ is compact w.r.t. the weak convergence of measures.

Given a lower semicontinuous cost function $k:\mathbb{H}^n\times\mathbb{H}^n \rightarrow [0,+\infty]$, the Monge-Kantorovich problem between $\mu$ and $\nu$, associated with the cost $k$,
\begin{equation} \label{e:MK} 
	\inf_{\gamma \in \Pi(\mu,\nu)} \int_{\mathbb{H}^n\times\mathbb{H}^n} k(x,y)\,d\gamma(x,y),
\end{equation}
admits solutions. We call 
these solutions \textit{optimal transport plans} for the cost function $k$ and we denote by
\begin{equation*}
	\Pi_k(\mu,\nu):=\{\gamma\in\Pi(\mu,\nu):\gamma \text{ solves } \eqref{e:MK}\}. 
\end{equation*}
This set is closed in $\Pi(\mu,\nu)$, w.r.t. the weak convergence of measures. Moreover, if $\gamma\in\Pi_k(\mu,\nu)$ and $\int_{\mathbb{H}^n\times\mathbb{H}^n}kd\gamma<+\infty$, then $\gamma$ is concentrated on a $k$-cyclically monotone set $\Gamma\subseteq\mathbb{H}^n\times\mathbb{H}^n$, i.e.

\begin{equation}\label{c-CM}
\sum_{i=1}^N k(x_i,y_i) \leq \sum_{i=1}^N k(x_{i},y_{\sigma(i)}),
\end{equation}
for any $N\geq 2$, any $(x_1,y_1), \dots, (x_N,y_N)\in\Gamma$ and any permutation $\sigma$ of $N$ elements. See  \cite[Chapter 1]{Santambrogiolibro}.

We say that a transport plan $\gamma\in \Pi(\mu,\nu)$ is \textit{induced by a map} if there exists a Borel map $T:\mathbb{H}^n\rightarrow \mathbb{H}^n$ such that $(I \otimes T)_\sharp\mu = \gamma$, where $(I \otimes T)(x) := (x,T(x))$. We will refer to a Borel map $T:\mathbb{H}^n\rightarrow \mathbb{H}^n$ solving \begin{equation}\label{mongepb}
	\inf_{T_{\#}\mu=\nu}\int_{\mathbb{H}^n}k(x,T(x))d\mu(x),
\end{equation}
as an \textit{optimal transport map} for the cost function $k$. If $\gamma\in\Pi_k(\mu,\nu)$ is induced by a map $T$, then $T$ is an optimal transport map for the cost $k$. Moreover, if any $\gamma\in\Pi_k(\mu,\nu)$ is induced by a map, then there exists a unique optimal transport map. Hence also $\gamma\in\Pi_k(\mu,\nu)$ is unique.



When $k=d_{CC}$, from the arguments above it follows that the Monge-Kantorovich problem
\begin{equation}\label{MKH}
	\inf_{\gamma\in \Pi(\mu,\nu)} \int_{\mathbb{H}^n\times\mathbb{H}^n} d_{CC}(x,y)\,d\gamma(x,y),
\end{equation}
admits at least a solution, concentrated on a $d_{CC}$-cyclically monotone set. From now on we will denote by
\begin{equation}\label{23luglio}
    \Pi_1(\mu,\nu):=\left\{\gamma\in\Pi(\mu,\nu):\gamma \textnormal{ solves }\eqref{MKH}\right\}.
\end{equation}

The explicit representation of geodesics, together with the absolute continuity of either $\mu$, or $\nu$, with respect to the Haar measure of the group, imply that any optimal transport plan (for the cost function $d_{CC}$) is concentrated on the set $E\subseteq\mathbb{H}^n\times\mathbb{H}^n$ of pairs of points connected by a unique geodesic, see \eqref{KAPPA} for its definition.
\begin{prop} \label{pi1.1}
    If either $\mu\ll\mathcal{L}^{2n+1}$, or $\nu\ll\mathcal{L}^{2n+1}$, and $\gamma\in \Pi_1(\mu,\nu)$, then for $\gamma$-a.e. $(x,y)\in\mathbb{H}^n\times\mathbb{H}^n$, there exists a unique geodesic between $x$ and $y$, i.e. $\gamma(\mathbb{H}^n\times\mathbb{H}^n\setminus E)=0$.
\end{prop}
See \cite[Lemma 4.1]{DePascale2} for the proof.

We denote by
\begin{equation}\label{lipset}
    \text{Lip}_1(\mathbb{H}^n, d_{CC}):=\left\{u:\mathbb{H}^n\to\mathbb{R}:|u(x)-u(y)|\leq d_{CC}(x,y),\forall x,y\in\mathbb{H}^n\right\}.
\end{equation}
The following important theorem holds.
\begin{teo}\label{1lip_potential}
There exists a function $u\in \textnormal{Lip}_1(\mathbb{H}^n,d_{CC})$ so that
\begin{equation*}
	 \min_{\gamma\in\Pi(\mu,\nu)}\int_{\mathbb{H}^n\times\mathbb{H}^n} d_{CC}(x,y)\,d\gamma(x,y) 
	= \int_{\mathbb{H}^n} u(x)\,d\mu(x) - \int_{\mathbb{H}^n} u(y)\,d\nu(y),
\end{equation*}
and 
$\gamma\in \Pi(\mu,\nu)$ is optimal if and only if 
\begin{equation*}
    u(x) - u(y) = d_{CC}(x,y),  \quad\text{for }\gamma-\text{a.e. }(x,y)\in \mathbb{H}^n\times\mathbb{H}^n.
\end{equation*}
\end{teo}

We call such a $u\in\text{Lip}_1(\mathbb{H}^n,d_{CC})$ a \textit{Kantorovich potential}.

From now on, we fix a Kantorovic potential $u\in\text{Lip}_1(\mathbb{H}^n,d_{CC})$ and we use it to check the optimality of transport plans. In this way one can select some optimal transport plans that satisfy a monotonicity condition, in the following sense.

If $\gamma\in\Pi_1(\mu,\nu)$, Theorem \ref{1lip_potential} and the Lipschitzianity of $u$ imply that 
\begin{equation*}\label{condition1}
	u(\sigma_{x,y}(t))=u(x)-d_{CC}\left(x,\sigma_{x,y}(t)\right),\quad\forall t\in[0,1],
\end{equation*}
for $\gamma$-a.e $(x,y)\in\mathbb{H}^n\times\mathbb{H}^n$ and any $\sigma_{x,y}$ geodesic between $x$ and $y$.
In this way one can define an order relation on $\sigma_{x,y}$ in the following way: let $t_1,t_2\in[0,1]$, $x'=\sigma_{x,y}(t_1)$ and $x'':=\sigma_{x,y}(t_2)$, then
\begin{equation}\label{orderrelation}
	x'\leq x''\Longleftrightarrow u(x')\geq u(x'').
\end{equation}
We denote by 
$\Pi_2(\mu,\nu)$ the set of transport plans solving the secondary variational problem
\begin{equation}\label{secvarpb}
	\inf_{\gamma \in \Pi_1(\mu,\nu)} \int_{\mathbb{H}^n\times\mathbb{H}^n} d_{CC}(x,y)^2\,d\gamma(x,y).
\end{equation}
This problem admits solutions since the functional $\gamma\mapsto\int_{\mathbb{H}^n\times\mathbb{H}^n}d_{CC}(x,y)^2\,d\gamma(x,y)$  is continuous w.r.t. the weak convergence of measures and $\Pi_1(\mu,\nu)$ is compact w.r.t. the same convergence. Theorem \ref{1lip_potential} allows us to rephrase this problem as a classical Monge-Kantorovich problem \eqref{e:MK} with cost $k(x,y) = \beta(x,y)$, where 
\begin{equation*}
	\beta(x,y) = \begin{cases}
		d_{CC}(x,y)^2 \quad \text{if } u(x)-u(y) = d_{CC}(x,y),\\
		+\infty \qquad \phantom{\text{if}} \text{otherwise}.
	\end{cases}
\end{equation*}
Since $\beta$ is lower semicontinuous and $\int_{\mathbb{H}^n\times\mathbb{H}^n} \beta(x,y) \,d\gamma(x,y)<+\infty$ for all $\gamma \in \Pi_2(\mu,\nu)$, it follows that any $\gamma \in \Pi_2(\mu,\nu)\subset\Pi_1(\mu,\nu)$ is concentrated on a $\beta$-cyclically monotone set $\Gamma$, i.e. 
\begin{equation}\label{var1}
	u(x) - u(y) = d_{CC}(x,y), \quad \forall (x,y)\in \Gamma,
\end{equation}
and
\begin{equation}\label{var2}
	\beta(x,y) +\beta(x',y') \leq \beta(x,y')+\beta(x',y),\quad \forall (x,y),(x',y')\in \Gamma.
\end{equation}

Using the non-branching property of $(\mathbb{H}^n,d_{CC})$ one can prove that geodesics used by an optimal transport plan cannot bifurcate. Moreover, if a transport plan solves also \eqref{secvarpb} then \eqref{var1} and \eqref{var2} imply a one-dimensional monotonicity condition along geodesics. More precisely, the following result \cite[Lemma 4.2 and Lemma 4.3]{DePascale2} holds.
\begin{prop}\label{monotone}
	Let $\gamma\in \Pi_1(\mu,\nu)$. Then, $\gamma$ is concentrated on a set $\Gamma$ such that for all $(x,y)$, $(x',y')\in \Gamma$ such that $x\not=y$ and $x\not=x'$, if $x'$ lies on a geodesic between $x$ and $y$ then all points $x$, $x'$, $y$ and $y'$ lie on the same geodesic.
	
	Moreover if $\gamma\in \Pi_2(\mu,\nu)$, then the condition $x<x'$ implies $y\leq y'$.
\end{prop}

As far as we know, in the Heisenberg Group has not been proven that any $\gamma\in\Pi_2(\mu,\nu)$ is induced by a map, and hence $\gamma\in\Pi_2(\mu,\nu)$ is unique. See \cite[Theorem 28]{Feldman} or \cite[Theorem 3.18]{Santambrogiolibro} for the analogous result in the Riemannian setting. Anyway in \cite{DePascale2} the authors proved that some particular transport plans in $\Pi_2(\mu,\nu)$, more precisely the ones that can be selected through the variational approximation below, are induced by maps.
\\

Following \cite[Section 5]{DePascale2}, we recall the aforementioned variational approximation procedure. Let $K$ be a compact subset of $\mathbb{H}^n$ such that 
\begin{equation*}
    \text{supp}(\mu)\cup\text{supp}(\nu)\subset K,
\end{equation*}
and let us denote by 
\begin{equation*}
	\Pi:=\{\gamma\in\mathcal{P}(\mathbb{H}^n\times\mathbb{H}^n): (\pi_1)_{\#}\gamma=\mu,\  \text{supp}((\pi_2)_\#\gamma)\subset K\}.
\end{equation*}

For any $\varepsilon\in\mathbb{R}_+$, we can consider the family of minimization problems
\begin{equation}\tag{$P_{\varepsilon}$}\label{varapprox} 
	\min\{C_\varepsilon(\gamma):\gamma\in\Pi\},
\end{equation}
where
\begin{multline*}
	C_{\varepsilon}(\gamma):=\frac{1}{\varepsilon}\, W_1((\pi_2)_\sharp\gamma,\nu) + \int_{\mathbb{H}^n\times\mathbb{H}^n} d_{CC}(x,y)\,d\gamma(x,y) \\ +  \varepsilon \int_{\mathbb{H}^n\times\mathbb{H}^n} d_{CC}(x,y)^2\,d\gamma(x,y) + \varepsilon^{6n+8} \textnormal{card}{(\textnormal{supp}({(\pi_2)_\sharp\gamma}))},
\end{multline*}
where $W_1((\pi_2)_\sharp\gamma,\nu)$ denotes the \textit{1-Wasserstein distance} between the two measures $(\pi_2)_\sharp\gamma$ and $\nu$,
\begin{equation*}
    W_1((\pi_2)_\sharp\gamma,\nu):=\min\left\{\int_{\mathbb{H}^n\times\mathbb{H}^n}d_{CC}(x,y)d\gamma(x,y):\gamma\in\Pi((\pi_2)_\sharp\gamma,\nu)\right\}.
\end{equation*}

For any $\varepsilon>0$, the minimization problem \ref{varapprox} admits at least one solution with finite value. Moreover the following result holds.

\begin{lemma} \label{optpi2}
Let $(\varepsilon_k)_{k\in\mathbb{N}}\subset\mathbb{R}_+$ be such that $\varepsilon_k\underset{k\rightarrow\infty}{\longrightarrow}0$ and $\gamma_{\varepsilon_k}$ be a solution to $(P_{\varepsilon_k})$. If $\gamma_{\varepsilon_k}\rightharpoonup\gamma\in\mathcal{P}(\mathbb{H}^n\times\mathbb{H}^n)$, then $\nu_{\varepsilon_k}:=(\pi_2)_\sharp \gamma_{\varepsilon_k}\rightharpoonup\nu$ and $\gamma \in \Pi_2(\mu,\nu)$.
\end{lemma}

The following theorem (\cite[Theorem 8.1]{DePascale2}) guarantees the existence of optimal transport plans induced by maps, and hence the existence of solutions to \eqref{mongepb}. In particular the optimal transport plans that turn out to be induced by maps are the weak limit of solutions $(\gamma_{\varepsilon_k})_{k\in\mathbb{N}}$ to $(P_{\varepsilon_k})_{k\in\mathbb{N}}$, for some $\varepsilon_k\to0$ as $k\to\infty$. Therefore they are monotone on the transport rays in the sense of \eqref{orderrelation}.

\begin{teo}\label{mainthmbis} If $\mu\ll\mathcal{L}^{2n+1}$, then there exists an optimal transport map $T:\mathbb{H}^n\rightarrow\mathbb{H}^n$ such that $\gamma=(\textnormal{Id}\otimes T)_\#\mu\in\Pi_2(\mu,\nu)$.  
\end{teo}

\section{Horizontal transport density on $\mathbb{H}^n$}\label{sectrandens}
In this section we introduce the notion of horizontal transport density, extending to the Heisenberg group the presentation provided in \cite{Santambrogio2}.  A horizontal transport density is a measure representing the amount of transport taking place along geodesics in each region of $\mathbb{H}^n$. In particular, we study conditions under which transport densities are Lebesgue absolutely continuous w.r.t. the Haar measure of the group, with $L^p$ density.

Let $\mu,\nu\in\mathcal{P}_c(\mathbb{H}^n)$ and let us fix a selection of geodesics
\begin{equation}\label{31luglio}
    S:\mathbb{H}^n\times\mathbb{H}^n \rightarrow \textnormal{Geo}(\mathbb{H}^n), 
\end{equation}
$S(x,y)=\sigma_{x,y}\in\textnormal{Geo}(\mathbb{H}^n)$, that is $\gamma$-measurable for any $\gamma\in\Pi_1(\mu,\nu)$, according to \eqref{19marzo1}. One can associate with any optimal transport plan $\gamma\in\Pi_1(\mu,\nu)$ a positive and finite Radon measure $a_\gamma\in\mathcal{M}_+(\mathbb{H}^n)$, defined as 
\begin{equation}\label{transport density}
	\int_{\mathbb{H}^n} \phi(x)da_\gamma(x):=\int_{\mathbb{H}^n\times\mathbb{H}^n}L_\phi(\sigma_{x,y})d\gamma(x,y),\quad\forall\phi\in C_c(\mathbb{H}^n,\mathbb{R}_+).
\end{equation}
Here $L_\phi(\sigma_{x,y})$ denotes the horizontal length of $\sigma_{x,y}$, weighted by $\phi$, see \eqref{3agosto} for its definition. The total mass of $a_\gamma$ satisfies
\begin{equation*}
    a_\gamma(\mathbb{H}^n)\leq\min_{\tilde\gamma\in\Pi(\mu,\nu)}\int_{\mathbb{H}^n\times\mathbb{H}^n}d_{CC}(x,y)d\tilde\gamma(x,y).
\end{equation*}

Moreover, the measure $a_\gamma$ is a compactly supported measure, see \cite{circelli2024continuous}. 
This measure is generally called \textit{horizontal transport density}.  By definition, it also follows that if $A$ is a Borel set, then
\begin{equation}\label{densitylength}
	a_\gamma(A)=\int_{\mathbb{H}^n\times\mathbb{H}^n}\mathcal{H}^1(A\cap S(x,y))d\gamma(x,y).
\end{equation}

Let us remark that if either $\mu\ll\mathcal{L}^{2n+1}$, or $\nu\ll\mathcal{L}^{2n+1}$, then $a_\gamma$ does not depend on the fixed selection $S$. See Proposition \ref{pi1.1}.

\subsection{Absolute continuity of horizontal transport densities}
The first goal is to prove the existence of at least one horizontal transport density that is absolutely continuous w.r.t. the Haar measure of the group.

Given $\gamma\in\Pi_1(\mu,\nu)$, we denote by $\mu_t$ the displacement interpolation between $\mu$ and $\nu$  
\begin{equation*}\label{interpolation}
	\mu_t:=((S_t)_\# \gamma)_{t\in [0,1]}.
\end{equation*}
Hence, the horizontal transport density $a_\gamma$ may be written as 
\begin{equation*}\label{density2}
	a_\gamma=\int_0^1(S_t)_\# (d_{CC}\ \gamma)dt,
\end{equation*}
where $d_{CC}\ \gamma$ is a positive Borel measure on $\mathbb{H}^n\times\mathbb{H}^n$.
Since $\mu$ and $\nu$ have bounded support, then there exists $C>0$ such that $d_{CC}(x,y)\leq C$, for any $(x,y)\in\textnormal{supp}(\gamma)$ and hence 
\begin{equation}
	\label{density3}
	a_\gamma\leq C\int_0^1\mu_t dt.
\end{equation}
In order to prove that $a_\gamma$ is absolutely continuous w.r.t. $\mathcal{L}^{2n+1}$, it is enough to prove that $\mu_t$ is absolutely continuous w.r.t. $\mathcal{L}^{2n+1}$, for almost every $t\in[0,1]$. In this way we would get that, whenever $\mathcal{L}^{2n+1}(A)=0$, then 
\begin{equation}\label{density4}
a_\gamma(A)\leq C\int_0^1\mu_t(A)dt=0.
\end{equation}
\\

We introduce now  the following lemma which guarantees that minimizing geodesics arising in 
optimal plans cannot intersect at intermediate points. This result will be useful in the proof of Theorem \ref{absolutecont}.

\begin{lemma} \label{pi1.3}
	Let $\gamma\in \Pi_1(\mu,\nu)$. Then $\gamma$ is concentrated on a set $\Gamma$ such that $\forall(x,y),(x',y')\in \Gamma$ with $(x,y)\not=(x',y')$, if two transport rays between these two pairs of points intersect at an interior point $z\in\mathbb{H}^n$, then all points $x$, $x'$, $y$, $y'$ and $z$ lie on the same transport ray.
	Moreover if $\gamma\in\Pi_2(\mu,\nu)$, then either $x\leq x'\leq z\leq y\leq y'$ or $x'\leq x\leq z\leq y'\leq y$.
\end{lemma}
\begin{proof}
	We first recall that \eqref{c-CM} reads as
	\begin{equation}\label{aaa}
		d_{CC}(x,y)+d_{CC}(x',y')\leq d_{CC}(x,y')+d_{CC}(x',y),
	\end{equation}
	$\forall (x,y),(x',y')\in\Gamma$.
	Let $\sigma:[0,d_{CC}(x,y)]\rightarrow\mathbb{H}^n$ be a geodesic between $x$ and $y$, $\tilde\sigma:[0,d_{CC}(x',y')]\rightarrow\mathbb{H}^n$ a geodesic between $x'$ and $y'$, $z\in\sigma(0,d_{CC}(x,y))\cap\tilde\sigma(0,d_{CC}(x',y'))$, so $z=\sigma(d_{CC}(x,z))=\tilde\sigma(d_{CC}(x',z))$. We denote by $\alpha$ the curve between $x$ and $y'$ defined in the following way:
    \begin{displaymath}
		\alpha(t):=\begin{cases}
			\sigma\left(\frac{d_{CC}(x,z)}{d_{CC}(x',z)}t\right),\quad &\textnormal{if}\,\ t\in[0,d_{CC}(x',z)],\\
			\tilde\sigma(t),\quad &\textnormal{if}\,\ t\in(d_{CC}(x',z),d_{CC}(x',y')].
		\end{cases}
    \end{displaymath}
    We will prove that $\alpha$ is geodesic between $x$ and $y'$. Indeed, otherwise we would have 
	\begin{multline}\label{aab}
		d_{CC}(x,y')<\ell_H(\alpha)=\ell_H(\alpha_{|[0,d_{CC}(x',z)]})+\ell_H(\alpha_{|[d_{CC}(x',z),d_{CC}(x',y')]})\\=d_{CC}(x,z)+d_{CC}(z,y').
	\end{multline}
	Since $z$ lies on both the geodesic between $x$ and $y$ and the geodesic between $x'$ and $y'$, it follows that 
	\begin{equation}\label{aac}
		\begin{cases}
			d_{CC}(x,y)=d_{CC}(x,z)+d_{CC}(z,y);\\
			d_{CC}(x',y')=d_{CC}(x',z)+d_{CC}(z,y').
		\end{cases}	
	\end{equation}
	By replacing (\ref{aac}) in (\ref{aab}), we obtain:
	\begin{equation}\label{aad}
		d_{CC}(x,y')+d_{CC}(z,y)+d_{CC}(x',z)<d_{CC}(x,y)+d_{CC}(x',y').
	\end{equation}
	By the triangle inequality follows that:
	\begin{displaymath}
		d_{CC}(x',y)\leq d_{CC}(x',z)+d_{CC}(z,y),
	\end{displaymath}
	and then, by replacing this last inequality in (\ref{aad}), we obtain
	\begin{displaymath}
		d_{CC}(x,y')+d_{CC}(x',y)<d_{CC}(x,y)+d_{CC}(x',y'),
	\end{displaymath}
	and this contradicts (\ref{aaa}). It follows that $\tilde\sigma$ and $\alpha$ are geodesics that coincide on the non-trivial interval $[d_{CC}(x',z),d_{CC}(x',y')]$. Since $\mathbb{H}^n$ is non-branching, this implies that $\tilde\sigma$ and $\alpha$ are sub-arcs of the same geodesic, namely $\alpha$ if $d_{CC}(x',z)\leq d_{CC}(x,z)$ and $\tilde\sigma$ otherwise, on which all points $x,x',z,y'$ lie.
	
	The thesis follows from Proposition \ref{monotone}.
\end{proof}

Given a map $T:\mathbb{H}^n\rightarrow\mathbb{H}^n$, from now on we will denote by
\begin{equation*}
	T_t:=S_t\circ(\text{Id}\otimes T):\mathbb{H}^n\rightarrow\mathbb{H}^n,\quad \forall t\in[0,1]
\end{equation*}
where $T_t(x)$ is the point at distance $td_{CC}(x, T(x))$ from $x$ on the selected geodesic $S(x,T(x))$ between $x$ and $T(x)$. In particular if $\gamma\in\Pi_1(\mu,\nu)$ is induced by a transport map, i.e. is of the form $\gamma:=(\text{Id}\otimes T)_\#\mu\in\Pi_1(\mu,\nu)$, then
\begin{equation*}
    \mu_t={(T_{t})}_{\#}\mu. 
\end{equation*}

The previous lemma allows to prove the following result.

\begin{prop}\label{absolutecontinterp}
If $\mu\ll\mathcal{L}^{2n+1}$, then there exists an optimal transport plan $\gamma\in\Pi_2(\mu,\nu)$ such that the measure
\begin{equation}\label{absolutecontinter1}
	\mu_t:=(S_t)_\#\gamma\ll\mathcal{L}^{2n+1},\quad \forall t\in[0,1).
\end{equation}
\end{prop}

\begin{proof}	
First we suppose that $\nu$ is finitely atomic, with atoms $(y^i)_{i=1}^M$. Let $\gamma\in\Pi_2(\mu,\nu)\subset\Pi_1(\mu,\nu)$, as in Theorem \ref{mainthmbis}, which is monotone in the sense of \eqref{orderrelation} and induced by a transport map $T$. Let us denote by $\Gamma\subseteq\mathbb{H}^n\times\mathbb{H}^n$ the set $\gamma$ is concentrated on and \eqref{var1} and \eqref{var2} hold.
	
We denote by $\Omega_i:=T^{-1}(\{y^i\})\cap \pi_1(\Gamma)$: obviously these sets are mutually disjoint and $\mu(\Omega)=1$, where $\Omega:=\bigcup_{i=1}^M\Omega_i$.
	
Now we denote by $\Omega_i(t):=T_t(\Omega_i)$: if we fix $t\in[0,1)$, then $\Omega_i(t)\cap\Omega_j(t)=\emptyset$ for every $i,j=1,\ldots,M$. Indeed, if $\exists\ z\in\Omega_i(t)\cap\Omega_j(t)$ then $\exists\ x^i\in \Omega_i$ and $x^j\in\Omega_j$ such that $(x^i,y^i), (x^j,y^j)\in \Gamma, (x^i,y^i)\not=(x^j,y^j)$ and the geodesics between these two pairs of points intersect at $z$. Since $\gamma\in\Pi_2(\mu,\nu)$, by Theorem \ref{pi1.3} we can suppose that $x^i,y^i,x^j,y^j,z$ belong to the same unit-speed geodesic and $x^i\leq x^j\leq z\leq y^i\leq y^j$. In particular this means, on the one hand, that $td_{CC}(x^i,y^i)=d_{CC}(x_i,z)\geq d_{CC}(x^j,z)=td_{CC}(x^j,y^j)$, hence $d_{CC}(x^i,y^i)\geq d_{CC}(x^j,y^j)$. On the other hand $(1-t)d_{CC}(x^i,y^i)=d_{CC}(z,y_i)\leq d_{CC}(z,y^j)=(1-t)d_{CC}(x^j,y^j)$, hence $d_{CC}(x^i,y^i)\leq d_{CC}(x^j,y^j)$. It follows that $d_{CC}(x^i,y^i)= d_{CC}(x^j,y^j)$ and hence $d_{CC}(x^i,z)=d_{CC}(x^j,z)$ and $d_{CC}(z,y^i)=d_{CC}(z,y^j)$, which in turn implies that $x^j=x^i$ and $y^i=y^j$ and gives a contradiction. However it may happen that $x^i=y^i$ or $x^j=y^j$. Let us suppose that $x^i=y^i=z$: the same computation above implies that $d_{SR}(x^j,y^j)=0$, which in turns implies that $y^i=y^j$ and gives a contradiction.
	
Remember also that $\mu$ is absolutely continuous and hence there exists a correspondence $\varepsilon\mapsto\delta=\delta(\varepsilon)$ such that 
\begin{equation*}
	\mathcal{L}^{2n+1}(A)<\delta(\varepsilon)\Rightarrow\mu(A)<\varepsilon.
\end{equation*}

Let $A\subset\mathbb{H}^n$ be a Borel set, $t\in[0,1)$, then $\mu_t:=(T_t)_{\#}\mu$ is concentrated on $T_t(\text{supp}(\mu))$ and
\begin{equation*}
	\mu_t(A)=\sum_{i=1}^{M}\mu_t(A\cap\Omega_i(t))=\sum_{i=1}^M\mu(T_t^{-1}(A\cap\Omega_i(t)))=\mu\left(\bigcup_{i=1}^M(T_t^{-1}(A\cap\Omega_i(t)))\right),
\end{equation*}
since the sets $T_t^{-1}(A\cap\Omega_i(t))\subseteq\Omega_i$ are disjoint.
We observe that for any $x\in\Omega_i$, $T_t(x)=S_t(x,y^i)$, hence by \eqref{MCP} follows that
\begin{equation*}
    \mathcal{L}^{2n+1}(U)\leq\frac{1}{(1-t)^{2n+3}}\mathcal{L}^{2n+1}(T_t(U)),
\end{equation*}
for any $U\subset\Omega_i$. This in turn implies that
\begin{equation*}
    \mathcal{L}^{2n+1}(T_t^{-1}(A\cap\Omega_i(t)))\leq\frac{1}{(1-t)^{2n+3}}\mathcal{L}^{2n+1}(A\cap \Omega_i(t)),
\end{equation*}
and so 
\begin{equation*}
    \mathcal{L}^{2n+1}\left(\bigcup_{i=1}^M(T_t^{-1}(A\cap\Omega_i(t)))\right)\leq\frac{1}{(1-t)^{2n+3}}\mathcal{L}^{2n+1}(A).
\end{equation*}
Hence, it is sufficient to suppose that $\mathcal{L}^{2n+1}(A)<(1-t)^{2n+3}\delta(\varepsilon)$ to get $\mu_t(A)<\varepsilon$.	This proves that $\mu_t\ll\mathcal{L}^{2n+1}$.
	
Now, if $\nu$ is not finitely atomic, we can take a sequence $(\nu_k)_{k\in\mathbb{N}}$ of atomic measures weakly converging to $\nu$, for instance as in Lemma \ref{optpi2}. For any $k\in\mathbb{N}$, we consider an optimal transport plan $\gamma_k\in\Pi_2(\mu,\nu_k)$ as in the first part of the proof. Hence, the sequence $(\gamma_k)_{k\in\mathbb{N}}$ weakly converges to some optimal transport plan $\gamma\in\Pi_2(\mu,\nu)$; moreover the sequence $(\mu_t^k)_{k\in\mathbb{N}}$ weakly converges to the corresponding $\mu_t:=(S_t)_{\#}\gamma$, thanks to Proposition \ref{pi1.1} and \cite[Lemma 7.3]{DePascale2}. Take a set $A$ such that $\mathcal{L}^{2n+1}(A)<(1-t)^{2n+3}\delta(\varepsilon)$. Since the Lebesgue measure is regular, $A$ is included in an open set $B$ such that $\mathcal{L}^{2n+1}(B)<(1-t)^{2n+3}\delta(\varepsilon)$. Hence $\mu_t^k(B)<\varepsilon,\forall k\in\mathbb{N}$. Passing to the limit and using Portmanteau's Theorem, see \cite[Theorem 2.1]{Billingsley}, we get 
$$\mu_t(A)\leq\mu_t(B)\leq\liminf_k \mu_t^k(B)\leq\varepsilon.$$
This proves that $\mu_t\ll\mathcal{L}^{2n+1}$.\end{proof}

Now we are able to find at least an optimal transport plan $\gamma\in\Pi_1(\mu,\nu)$ such that the interpolation measures $\mu_t$ constructed from $\gamma$ are absolutely continuous for $t<1$.
\\

\begin{teo}\label{absolutecont}
If $\mu\ll\mathcal{L}^{2n+1}$, then there exists an optimal transport plan $\gamma\in\Pi_2(\mu,\nu)$ such that the measure $a_{\gamma}\ll\mathcal{L}^{2n+1}$.
\end{teo}
\begin{proof}
Let $\gamma\in\Pi_2(\mu,\nu)$ satisfying \eqref{absolutecontinter1}. Then, the thesis follows immediately from \eqref{density4} applied to $a_{\gamma}$.
\end{proof}

Obviously the previous argument depends only on one of the two marginals and it is completely symmetric: if $\nu\ll\mathcal{L}^{2n+1}$, again one can get the existence of an optimal transport plan $\gamma\in\Pi_2(\mu,\nu)$ such that the associated horizontal transport density $a_{\gamma}$ is absolute continuous w.r.t. the $(2n+1)$-dimensional Lebesgue measure.


\subsection{$p$-summability of horizontal transport densities}

In this subsection we prove the existence of at least one horizontal transport density belonging to $L^p$, for some values of $p$.

From now on, given $\lambda\in\mathcal{M}_+(\mathbb{H}^n)$ we will write that $\lambda\in L^p$ if $\lambda\ll\mathcal{L}^{2n+1}$, with density $\rho\in L^p$. We will denote by $\|\lambda\|_p:=\|\rho\|_{L^p}$.

Let $\gamma\in\Pi_1(\mu,\nu)$ as in Theorem \ref{absolutecont}. From \eqref{density3} and the Minkowski inequality it follows that
\begin{equation}\label{Minkowski}
	\|a_\gamma\|_p\leq C\int_{0}^1\|\mu_t\|_pdt.
\end{equation}

In order to prove $p$-summability of $a_\gamma$, it is enough to estimate the $L^p$ norm of $\mu_t$ as a function of the variable $t$. This will be established in the following  theorem: 
\begin{prop}\label{summabilityinterp}
If $\mu\in L^p$, for some $p\in[1,\infty]$, then there exists an optimal transport plan $\gamma\in\Pi_2(\mu,\nu)$ such that $\mu_t:=(S_t)_{\#}\gamma\in L^p$ and 
\begin{equation}\label{bb3}
	\|\mu_t\|_p\leq (1-t)^{-(2n+3)/q}\|\mu\|_p,\quad \forall t\in[0,1),
\end{equation}
 where $q:=\frac{p}{p-1}$.
\end{prop}
\begin{proof}
Let us denote by $\rho$ the density of $\mu$ w.r.t. $\mathcal{L}^{2n+1}$. Let us consider first the discrete case: let us assume that the target measure $\nu$ is finitely atomic and let us denote by $(y^i)_{i=1,\ldots,M}$ its atoms. Let us consider an optimal transport plan $\gamma\in\Pi_2(\mu,\nu)$, as in the proof of Proposition \ref{absolutecontinterp}, concentrated on some set $\Gamma$. Since $\gamma$ is induced by a map $T$, we denote by $\Omega_i:=T^{-1}(\{y^i\})\cap\pi_1(\Gamma)$, for $i\in\{1,\ldots,M\}$, so that for $\gamma$-a.e. $(x,y)\in\Omega_i\times\mathbb{H}^n$, we have $y=y^i$. Let us consider the corresponding interpolation measures $\mu_t\ll\mathcal{L}^{2n+1}$ for every $t\in[0,1)$; moreover, for all $\phi\in C_c(\mathbb{H}^n,\mathbb{R}_+)$, by definition of push-forward we get that 
\begin{align*}
    \int\phi(x)d\mu_t(x)=&\sum_{i=1}^M\int_{\Omega_i}\phi(S_t(x,y_i))d\gamma(x,y_i)=\\
	=&\sum_{i=1}^M\int_{\Omega_i}\phi(T_t(x))d\mu(x).
\end{align*}
Let us fix $i\in\{1,\ldots,M\}$ and let us denote by $\rho_t$ the density of $\mu_t$ w.r.t. $\mathcal{L}^{2n+1}$ and by $\rho_t^i:={\rho_t}_{\lfloor \Omega_i}$. Let us take the change of variable $z=S_t(x,y_i)={T_t}_{\lfloor\Omega_i}(x)$. We know from Lemma \ref{pi1.3} and the disjointness of the sets $\Omega_i(t)$ that this map is injective. Then, for all $\phi\in C_c(\mathbb{H}^n,\mathbb{R}_+)$ we get 
\begin{align*}
	\int_{\Omega_i}\phi(x)d\mu^i_t(x)&=\int_{\Omega_i}\phi(T_t(x))\rho(x)dx=\\&=
		\int_{\Omega_i(t)}\phi(z)\rho(T_t^{-1}(z))|\det D_x(S_t(x,y^i))|^{-1}dz.
\end{align*}
Hence, we have that 
\begin{equation*}
    \rho_t^i(z)=\rho(T_t^{-1}(z))|\det D_x(S_t(x,y^i))|^{-1},\quad \text{for a.e. } z\in\Omega_i(t).
\end{equation*}
Consequently, we get 
\begin{align*}
	\|\rho_t^i\|^p_{L^p(\Omega_i(t))}&=\int_{\Omega_i(t)}\rho(T_t^{-1}(z))^p|\det D_x(S_t(x,y^i))|^{-p}dz=\\&=\int_{\Omega_i}\rho(x)^p|\det D_x(S_t(x,y^i))|^{1-p}dx.
\end{align*}
Hence from \eqref{det} it follows that
\begin{align*}
    \|\rho_t^i\|^p_{L^p(\Omega_i(t))}\leq (1-t)^{(1-p)(2n+3)}\|\rho\|^p_{L^p(\Omega_i)},\quad \forall i\in\{1,\ldots,M\}.
\end{align*}
Then, we have
\begin{equation}\label{discretepsumm}
	\|\mu_t\|_p\leq (1-t)^{-(2n+3)/q}\|\mu\|_p,\quad \forall t\in(0,1),
\end{equation}
where $q:=\frac{p}{p-1}$.

If $\nu$ is not finitely atomic, again we take a sequence $(\nu_k)_{k\in\mathbb{N}}$ of atomic measures weakly converging to $\nu$, for instance as in Lemma \ref{optpi2}. We consider a sequence $(\gamma_k)_{k\in\mathbb{N}}\subset\Pi_2(\mu,\nu_k)$ of optimal plans satisfying \eqref{discretepsumm}: this sequence weakly converges to an optimal plan $\gamma\in\Pi_2(\mu,\nu)$ and $\mu_t^k$ weakly converge to the corresponding $\mu_t:=(S_t)_{\#}{\gamma}$, see again Proposition \ref{pi1.1} and \cite[Lemma 7.3]{DePascale2}. Hence, we get that
\begin{equation*}
    \|\mu_t\|_p\leq\liminf_{k \rightarrow 0}\|\mu_t^k\|_p\leq (1-t)^{-(2n+3)/q}\|\mu\|_p.
\end{equation*}
\end{proof}

Now we are able to prove the following theorem.
\begin{prop}\label{summability1}
If $\mu\in L^p$, for some $p\in[1,\infty]$, the following results hold: if $p<\frac{2n+3}{2n+2}$, then there exists $\gamma\in\Pi_2(\mu,\nu)$ such that $a_{\gamma}\in L^p$; otherwise, there exists $\gamma\in\Pi_2(\mu,\nu)$ such that $a_{\gamma}\in L^s$, for $s<\frac{2n+3}{2n+2}$.

\end{prop}

\begin{proof}
Let $\gamma\in\Pi_2(\mu,\nu)$ satisfying \eqref{bb3}. Then, it follows from \eqref{Minkowski} applied to $a_{\gamma}$ that
\begin{equation*}
	\|a_{\gamma}\|_p\leq C\int_0^1 \|\mu_t\|_pdt\leq C\|\mu\|_p\int_0^1(1-t)^{-(2n+3)/q}dt.
\end{equation*}
The last integral is finite whenever $q>2n+3$, i.e. $p<\frac{2n+3}{2n+2}$.
	
If $p\geq\frac{2n+3}{2n+2}$ the thesis follows from the fact that any density in $L^p$ also belongs to any $L^s$ space for $s<p$.
\end{proof}

If also $\nu\in L^p$ then, by symmetry, one can find an optimal transport plan $\tilde\gamma\in\Pi_2(\mu,\nu)$, possibly different from the one in Proposition \ref{summabilityinterp}, such that $\tilde{\mu}_t:=\left(S_t\right)_\#\tilde\gamma\in L^p$ and it satisfies
\begin{equation}\label{bb1}
	\|\tilde{\mu}_t\|_p\leq t^{-(2n+3)/q}\|\nu\|_p,\quad \forall t\in(0,1].
\end{equation}

In the Euclidean setting  
 \cite[Theorem 3.18]{Santambrogiolibro} or \cite{Feldman2}, and in the more general Riemannian setting, see \cite{Feldman}, $\Pi_2(\mu,\nu)$ consists of a unique element, so that it is possible to glue together \eqref{bb3}, for $t\leq\frac{1}{2}$, and \eqref{bb1}, for $t\geq\frac{1}{2}$, and get the existence of a horizontal transport density in $L^p$. Unfortunately, this uniqueness result is still an open problem in the Heisenberg group, hence we cannot glue together \eqref{bb3} and \eqref{bb1}, and deduce anything about the summability of $a_\gamma$.

\section{Weighted distance induced by a $L^q$ density}
A measure with continuous density naturally induces a weighted length of curves, leading to the definition of weighted distance. The aim of this section is to extend the definition of weighted distance to measures with $L^q$ densities. While the statements of the results appear similar to those in the Euclidean setting (contained in \cite{Santambrogio1}), their proofs differ significantly due to the geometric properties of the space.

We consider a bounded domain $\Omega\subset\mathbb{H}^{n}$ and we assume that its boundary has $C^{1,1}$ regularity, in the Euclidean sense. We consider the space 
\begin{equation}\label{acca}
	H:=\big{\{}\sigma\in \textnormal{Lip}([0,1],\overline{\Omega}):\ \sigma\ \textnormal{is\ horizontal}\big{\}}	 
\end{equation}
of horizontal curves on $\overline{\Omega}$ parameterized on $[0,1]$, viewed as subset of $C([0,1],\overline{\Omega})$ equipped with the topology of uniform convergence. Given $x,y\in\overline\Omega$, we denote by 
\begin{equation}\label{horcrvxy}
    H^{x, y}:=\{\sigma\in H:\ \sigma(0)=x,\ \sigma(1)=y\}.
\end{equation}
The regularity assumption on the boundary of $\Omega$ guarantees the existence of at least one horizontal curve between any pair of points in $\overline\Omega$. This assumption replaces the stronger convexity assumption used in the Euclidean setting in \cite{Santambrogio1},  since non-trivial geodesic convex sets do not exist in the Heisenberg group, as proved in \cite{Monti2}. Therefore, the following important lemma holds.

\begin{lemma}\label{nonemptinesshorcur}
The set $H^{x,y}\not=\emptyset$, for any $x,y\in\overline{\Omega}$. 
\end{lemma}

\begin{proof}
Let $x,y\in\Omega$. Since $\Omega$ is open and connected,  then the statement follows from the Chow-Rashevsky Theorem applied to the manifold $\Omega$, equipped with
the sub-Riemannian structure inherited from $\mathbb{H}^n$, see \cite[Theorem 3.31]{barilariagrachev}.

Now, let us denote by $\mathcal{C}(\partial\Omega)$ the set of \textit{characteristic points} of $\partial\Omega$, that is 
\begin{equation*}    
\mathcal{C}(\partial\Omega):=\left\{x\in\partial\Omega:T_x\partial\Omega=H_x\mathbb{H}^n\right\},
\end{equation*}
and let us remark that at a non-characteristic point $z$ the fiber $H_z\mathbb{H}^n$  can be represented as 
$$H_z\mathbb{H}^n=H_z\partial\Omega \oplus \mathrm{span}\{\textbf{n}_H(z)\}.$$
Here, $H_z\partial\Omega = T_z\partial\Omega \cap H_z\mathbb{H}^n$ denotes the horizontal tangent space to $\partial\Omega$ at the point $z$, and $\textbf{n}_H(z)$ is the horizontal normal at the point $z\in\partial\Omega$, that is the
orthogonal projection of the Euclidean normal to $\Omega$ at $z$ on $H_z\mathbb{H}^n$. Then, for any horizontal vector field $Z$, the vector $Z(z)\in H_z\mathbb{H}^n$ admits an unique projection $\pi(Z(z))$ on the space $H_z\partial\Omega$. One can consider the horizontal vector field $\pi(Z):\Omega\to H\Omega$, $$z\mapsto \pi(Z)(z):=\pi(Z(z)),$$ 
where $\pi(Z)(z)=Z(z)$ if $z$ is a characteristic point.

Let us consider the case $x\in\partial\Omega$ and $y\in\Omega$. If $x\in\partial\Omega$ is a non-characteristic point, then the horizontal normal $\textbf{n}_H$ at the point $x$ does not vanish. As a result, if $\textbf{n}_H(x):=\sum_{i=1}^{2n}\textbf{n}_iX_i(x)$, one can consider the horizontal vector field $Z:=\sum_{i=1}^{2n}\textbf{n}_iX_i\in\mathfrak{h}_1^n$, $\delta >0$ and the horizontal curve
\begin{equation}
    \sigma(t):=\exp(-tZ)(x), 
\end{equation}
such that $\sigma\left([0, \delta]\right)\subset\overline\Omega$ and $z:=\sigma(\delta)\in\Omega$. Now one can consider a horizontal curve between $z$ and $y$, fully contained in $\Omega$.

If $x\in\mathcal{C}(\partial\Omega)$ is a characteristic point, from \cite[Theorem 1.2]{balogh2003size} it follows that the Hausdorff dimension w.r.t. the Euclidean metric is $\dim_E \mathcal{C}(\partial\Omega)<2n$. Then, there exists $v=\sum_{i=1}^{2n}v_iX_i(x)\in T_x\partial\Omega=H_x\mathbb{H}^n$ and $\delta>0$ such that the horizontal curve 
\begin{equation}
    \sigma(t)=\exp\Big(t\sum_{i=1}^{2n}v_i\pi(X_i)\Big)(x)  
\end{equation}
is well-defined, it belongs to $\partial\Omega$ and it is non-characteristic for all $t\in(0, \delta]$. Hence, one can find a horizontal curve between $z=\sigma(\delta)$ and $y\in\Omega$, using the argument above.

If $x,y\in\partial\Omega$, using the previous arguments one can connect them with $x',y'\in \Omega$ and then find a horizontal curve between $x'$ and $y'$, which is contained in $\Omega$.
\end{proof}

Given $\sigma\in H$ and $\phi\in C(\overline{\Omega},\mathbb{R}_+)$, we 
recall that $L_{\phi}(\sigma)$ is the horizontal length of the curve $\sigma$, weighted by $\phi$, introduced in \eqref{3agosto}.

Following \cite[Lemma 2.7]{Santambrogio1} one can prove that, for any $\phi \in C(\overline{\Omega},\mathbb{R}_+)$ and for any $\sigma\in H$, it holds that 
\begin{multline*}\label{LPHI}
	L_{\phi}(\sigma):=\sup\bigg{\{}\sum_{i=1}^n \left(\inf_{[t_i,t_{i+1}]}(\phi\circ\sigma)\right)d_{CC}(\sigma(t_i), \sigma(t_{i+1})):([t_i, t_{i+1}])_i \mbox{ is a partition of }[0,1]\bigg{\}}. 
\end{multline*}
Hence, the functional $H\ni\sigma\mapsto L_{\phi}(\sigma)$ is lower semicontinuous, and hence Borel, w.r.t. the uniform convergence. 

Given a weight $\phi\in C(\overline{\Omega},\mathbb{R}_+)$, we denote its associated weighted distance by
\begin{equation}\label{ccontinuous}
	c_\phi(x,y):=\inf\{L_{\phi}(\sigma)\,:\,\sigma\in H^{x,y}\},\quad \forall (x,y)\in\Omega\times\Omega.
\end{equation}   

Let us explicitly recall that the cost function $c_\phi(x,y)$
is a distance, if $\phi$ is continuous and  strictly positive, and it is only a pseudo distance, if $\phi$ is non-negative. 
In order to extend this pseudo distance to summable functions  $\phi$, we start with an estimate of the regularity of $c_\phi(x,y)$ in terms of the $L^q$ norm of $\phi$, 
where $q:=\frac{p}{p-1}$ is the dual exponent of some $p\in(1,+\infty)$. The proof is inspired by \cite[Proposition 3.2]{Santambrogio1} but requires many non trivial changes, due to the geometric structure of $\mathbb{H}^n$. 

\begin{prop}\label{cxicomp1}
If $q>N$, then there exists $C>0$ such that for any $\phi\in C(\overline\Omega,\mathbb{R}_+)$ and any $(x,y), (x',y')\in\Omega\times\Omega$, it holds
\begin{equation}\label{holderest}
		\vert c_{\phi}(x,y)-c_{\phi}(x',y')\vert \leq C\Vert \phi\Vert_{L^{q}(\Omega)} \left(d_{CC}(x,x')^{\alpha}+ d_{CC}(y,y')^{\alpha} \right),
	\end{equation}
where $\alpha:=1-\frac{N}{q}$.

\end{prop}
\begin{proof}
   
Let $\phi\in C(\overline{\Omega},\mathbb{R}_+)$ and $x,y\in\Omega$. For $k>0$ let $\sigma_k\in H^{x,y}$ be such that
\begin{equation*}
	\int_0^1 \phi(\sigma_k(t))\vert \dot \sigma_k(t)\vert_H dt\leq c_{\phi}(x,y)+\frac{1}{k}.
\end{equation*}
 In order to study the regularity of $c_\phi$ with respect to the second variable $y$, we choose a point $z_\varepsilon$ 
 which can be connected to $y$ by an horizontal segment: i.e. we fix a horizontal vector field $Z\in\mathfrak{h}^{n}_1$, such that $|Z|_H=1$, and we choose for all $\varepsilon>0$ the points and we choose for all $\varepsilon>0$ the points 
\begin{equation*}
	z_\varepsilon:=\exp\left(\varepsilon Z\right)(y),
\end{equation*}
such that $z_\varepsilon\in \Omega$. Now we modify the curve $\sigma_k$ into a curve  $\sigma_{k,t_0}\in H^{x,z_\varepsilon}$: we choose $t_0\in (0,1)$ and define
\begin{equation*}
    \sigma_{k,t_0}(t):=
    \begin{cases}
		\sigma_k\big{(}{\frac{t}{t_0}}\big{)} &\mbox{ if }t\in[0,t_0]\\  
		\tilde{\sigma}_{\varepsilon,y}\big{(}\frac{t-t_0}{1-t_0}\big{)} &\mbox{ if }t\in]t_0,1],
    \end{cases}
\end{equation*}
where
\begin{equation*}
    \tilde{\sigma}_{\varepsilon, y}(t)=\exp\left(t(\varepsilon Z)\right)(y),\quad t\in[0,1].
\end{equation*}
We then have, for all $k>0$
\begin{align*}
    c_\phi(x,z_\varepsilon)&\leq \int_0^1\phi(\sigma_{k,t_0}(t))\vert\dot\sigma_{k,t_0}(t)\vert_H dt=\\
	&=\int_0^1\phi(\sigma_k(t))\vert \dot\sigma_k(t)\vert_H dt+\int_0^1\phi(\tilde{\sigma}_{\varepsilon, y}(t))\vert \dot{\tilde{\sigma}}_{\varepsilon, y}(t)\vert_H dt\leq\\
	&\leq c_\phi(x,y)+\frac{1}{k}+\varepsilon\int_0^1 \phi(\tilde{\sigma}_{\varepsilon, y}(t))dt.
\end{align*}
Now, if $k\rightarrow+\infty$, we get
\begin{equation*}
    \frac{1}{\varepsilon}\left[c_{\phi}\left(x,\exp\left(\varepsilon Z\right)(y)\right)-c_{\phi}(x,y)\right]\leq \int_0^1 \phi(\tilde{\sigma}_{\varepsilon, y}(t))dt,
\end{equation*}
and, by similar argument:
\begin{equation*}
    \frac{1}{\varepsilon}\left[c_{\phi}(x,y)-c_{\phi}\left(x,\exp\left(\varepsilon Z\right)(y)\right)\right]\leq \int_0^1 \phi(\tilde{\sigma}_{\varepsilon, y}(1-t))dt.
\end{equation*}

Integrating with respect to $y$, raising to the power $q$ and using the fact that the function $y\mapsto\tilde{\sigma}_{\varepsilon, y}(t) $ has Jacobian determinant $1$, 
we get that, for any fixed $x$, $Zc_{\phi}(x,\cdot)\in L^{q}(\Omega)$ and $\|Zc_\phi(x,\cdot)\|_{L^q(\Omega)}\leq C\|\phi\|_{L^q(\Omega)}$. Since this holds for every $Z$, we have $c_{\phi}(x,\cdot)\in HW^{1,q}(\Omega)$, see \eqref{horsob}, and
\begin{equation}\label{EST}
    \|\nabla_Hc_{\phi}(x,\cdot)\|_{L^q(\Omega)}\leq\|\phi\|_{L^q(\Omega)},\quad \forall x\in\Omega.
\end{equation}
By symmetry we also get that
\begin{equation}\label{ESTx}
    \|\nabla_Hc_{\phi}(\cdot,y)\|_{L^q(\Omega)}\leq\|\phi\|_{L^q(\Omega)},\quad \forall y\in\Omega.	
\end{equation}

Since $q>N$ then if follows from \eqref{EST}, \eqref{ESTx} and \cite[Theorem 1.11]{garofalo1996isoperimetric}, that there exists $C>0$ such that
\begin{align*}
    \vert c_{\phi}(x,y)-c_{\phi}(x,y')\vert \leq C\Vert \phi\Vert_{L^{q}(\Omega)}  d_{CC}(y,y')^{\alpha},\quad  \forall x, y, y'\in\Omega,\\
    \vert c_{\phi}(x,y)-c_{\phi}(x',y)\vert \leq C\Vert \phi\Vert_{L^{q}(\Omega)}  d_{CC}(x,x')^{\alpha},\quad   \forall x, x', y\in\Omega,	
\end{align*}
with $\alpha=1-\frac{N}{q}.$ This proves \eqref{holderest}. 
\end{proof}

\begin{prop}\label{12aprile}
    If $q>N$, then for any $\phi\in C(\overline\Omega, \mathbb{R}_+)$, the function $c_\phi$ defined in \eqref{ccontinuous} admits a unique continuous extension as a function
\begin{equation*}
    c_\phi:\overline\Omega\times\overline\Omega\to\mathbb{R}_+,
\end{equation*}
with the same modulus of continuity. Moreover the definition at \eqref{ccontinuous} extends to all pairs $(x,y)\in\overline{\Omega}\times\overline{\Omega}$.
\end{prop}

\begin{proof}
The first part of the proof easily follows from \eqref{holderest}.

Let now consider $\phi\in C(\overline\Omega,\mathbb{R}_+)$, $\varepsilon_0>0$ and the continuous function $\phi+\varepsilon_0>0$. Given $(x,y)\in\overline\Omega\times\overline\Omega$ and $\left(x_n,y_n\right)_{n\in\mathbb{N}}\subset\Omega\times\Omega$, $\left(x_n,y_n\right)\to(x,y)$, we have   
\begin{equation*}
    c_{\phi+\varepsilon_0}(x,y):=\lim_{n\to+\infty}c_{\phi+\varepsilon_0}(x_n,y_n).
\end{equation*}
It means that $\forall \varepsilon>0$, there exists $n=n(\varepsilon)$ such that
\begin{equation*}
    \left|c_{\phi+\varepsilon_0}(x,y)-c_{\phi+\varepsilon_0}(x_n,y_n)\right|<\varepsilon,\quad\forall n>n(\varepsilon).
\end{equation*}
By definition of $c_{\phi+\varepsilon_0}(x_n,y_n)$ and by the invariance of $L_{\phi+\varepsilon_0}$ under reparametrization, there exists $\sigma_n\in\tilde{H}^{x_n,y_n}$ such that
\begin{equation*}
    |c_{\phi+\varepsilon_0}(x_n,y_n)-L_{\phi+\varepsilon_0}(\sigma_n)|<\varepsilon.
\end{equation*}
Hence, for any $n>n(\varepsilon)$ it holds that
\begin{equation*}
    \varepsilon_0|\dot{\sigma}_n|_H=\varepsilon_0\ell_{H}(\sigma_n)\leq L_{\phi+\varepsilon_0}(\sigma_n)<c_{\phi+\varepsilon_0}(x_n,y_n)+\varepsilon\leq M+\varepsilon,
\end{equation*}
where $M\geq0$. Then, the Ascoli-Arzelà Theorem implies that $\left(\sigma_n\right)_{n>n(\varepsilon_0)}\subset H$ admits a subsequence $\sigma_{n_k}\to\sigma$ uniformly as $k\to+\infty$. From \cite[Theorem 3.41]{barilariagrachev} it follows that $\sigma\in H^{x,y}$ and the lower semicontinuity implies $L_{\phi+\varepsilon_0}(\sigma)\leq\liminf_{k\to+\infty}L_{\phi+\varepsilon_0}(\sigma_{n_k})=c_{\phi+\varepsilon_0}$. Then, we can conclude that $\forall\varepsilon_0>0$, $\forall (x,y)\in\overline\Omega\times\overline\Omega$, there exists $\sigma\in H^{x,y}$, such that 
\begin{equation*}
    L_{\phi+\varepsilon_0}(\sigma)\leq c_{\phi+\varepsilon_0}(x,y).
\end{equation*}
Moreover,
\begin{equation*}
    L_{\phi}(\sigma)+\varepsilon_0\ell_{H}(\sigma)=L_{\phi+\varepsilon_0}(\sigma)\leq c_{\phi+\varepsilon_0}(x,y)=\lim_{n\to+\infty}c_{\phi+\varepsilon_0}(x_n,y_n)\leq\lim_{n\to+\infty}c_{\phi}(x_n,y_n)+O(\varepsilon_0),
\end{equation*}
hence, letting $\varepsilon_0\to0$, $L_\phi(\sigma)\leq c_\phi(x,y)$.

It remains to prove that $c_{\phi}(x,y)\leq L_{\phi}(\sigma)$, for any $\sigma\in H^{x,y}$.

Let us suppose that $x\in\Omega$ and $y\in\partial\Omega$, all the other cases can be deduced from this one. Let us consider an arbitrary horizontal curve $\sigma\in H^{x,y}$ and, for any $n\in\mathbb{N}$, we take $y_n\in B\left(y,\frac{1}{n}\right)\cap\Omega$. From Remark \ref{nonemptinesshorcur} it follows that there exists a horizontal curve $\sigma_n\in H^{y,y_n}$, such that $\ell_{H}(\sigma_n)\leq 2d_{CC}(y,y_n)$. Hence, $L_\phi(\sigma_n)\leq 2\|\phi\|_\infty d_{CC}(y,y_n)$. Given $t_0\in[0,1]$, we denote by $\tilde\sigma_{n,t_0}\in H^{x,y_n}$ the horizontal curve defined as 
\begin{equation*}
    \tilde\sigma_{n,t_0}(t):=\begin{cases}
        \sigma\left(\frac{t}{t_0}\right),&\text{if }t\in[0,t_0],\\
        \sigma_n\left(\frac{t-t_0}{1-t_0}\right),&\text{if }t\in[t_0,1].
    \end{cases}
\end{equation*}

From the invariance of the weighted sub-Riemannian length under reparametrization, it follows that
\begin{equation*}
    c_\phi(x,y) := \lim_{n\to\infty} c_\phi(x,y_n)\leq \liminf_{n\to\infty} L_\phi (\tilde \sigma_{n,t_0})=\liminf_{n\to\infty}\left(L_\phi (\sigma) + L_\phi(\sigma_n)\right)= L_\phi(\sigma).
\end{equation*}
Since $\sigma\in H^{x,y}$ is arbitrary, the thesis follows.
\end{proof}

\begin{cor}\label{cxicomp}
Let $(\phi_n)_{n\in\mathbb{N}}\subset  C(\overline\Omega,\mathbb{R}_+)$ be a bounded sequence in $L^{q}$, then the sequence $(c_{\phi_n})_{n\in\mathbb{N}}$ admits a subsequence that converges in $C(\overline\Omega\times\overline\Omega)$.
\end{cor}

\begin{proof}
    The existence of a subsequence of $(c_{\phi_n})_{n\in\mathbb{N}}$ that converges in $C(\overline\Omega\times\overline\Omega)$ follows from Ascoli-Arzelà's theorem. Indeed, equicontinuity follows from Proposition \ref{12aprile}, while the pointwise boundness is a consequence of the identity $c_{\phi_n}(x,x)=0$ and \eqref{holderest}.
\end{proof}

From now on, we assume that $q>N$. The next goal is to give an equivalent definition for $c_\phi$, that extends the notion for positive functions in $L^{q}$.

\begin{prop}\label{ccoincid}
	If $\phi\in C(\overline\Omega,\mathbb{R}_+)$, then
\begin{equation*}
	c_\phi(x,y)=\sup\left\{c(x,y)\,:c\in\mathcal{C}(\phi)\right\},\quad \forall(x,y)\in\overline\Omega\times\overline\Omega,
\end{equation*}
where 
\begin{equation}\label{Setcphi}
\mathcal{C}(\phi)=\left\{c=\lim_{n\rightarrow+\infty} c_{\phi_n}\,\textnormal{ in }C(\overline{\Omega}\times\overline{\Omega})\,:\,(\phi_n)_{n\in\mathbb{N}}\subset C(\overline{\Omega},\mathbb{R}_+),\phi_n\to\phi\,\textnormal{ in }L^{q}\right\}.
\end{equation}
\end{prop}

We first state two technical remarks that will be useful in the proof.

\begin{Remark}
Note that if we have a constant coefficient unitary 
horizontal vector 
 $W_1:=a_1X_1+\ldots +a_nX_n+a_{n+1}X_{n+1}+\ldots+a_{2n}X_{2n}\in\mathfrak{h}_1^1$, 
 it is  possible to perform a change of variable which sends the 
 vector $W_1$ to the first element of the canonical orthonormal basis. 
 Indeed, if we  denote by $W_2,\ldots,W_{2n}$ a basis of orthogonal complement $W_1^\perp$ in $\mathfrak{h}_1^1$ with respect to $\left\langle\cdot,\cdot\right\rangle_H$, and by $x$ a point, we can consider the change of variables $\Psi:\mathbb{R}^{2n+1}\to\mathbb{H}^{n},$
\begin{equation}\label{changeofcoordinates}
    \Psi(e_1,\ldots,e_{2n+1})=\exp(e_1W_1)\exp\left(\sum_{i=2}^{2n} e_{i}W_{i}+e_{2n+1}X_{2n+1}\right)(x).
\end{equation}
The pullback of the vector field $W_1$ by $\Psi$ is $\Psi_*W_1=\partial_{e_1}$ and the point $x$ will be the origin in the new  coordinate system.
\end{Remark}

\begin{Remark}\label{30luglio}
Any function $c\in \mathcal{C}(\phi)$ is a pseudo distance. Indeed, let $(\phi_n)_{n\in\mathbb{N}}\subset C(\overline\Omega,\mathbb{R}_+)$, such that $\phi_n\to\phi$ in $L^q$ and $c_{\phi_n}\to c$ in $C(\overline\Omega\times\overline\Omega)$. We know that the function $c_{\phi_n}$ is a pseudo distance, for any $n\in\mathbb{N}$. Hence, the thesis follows passing to the limit. In particular, for any $x,y,z\in\overline\Omega$ it holds
\begin{equation}\label{ineq}
		c(x,z)=\lim_{n\to+\infty}c_{\phi_n}(x,z)\leq\lim_{n\to+\infty}\left(c_{\phi_n}(x,y)+c_{\phi_n}(y,z)\right)=c(x,y)+c(y,z).
	\end{equation}
\end{Remark}

\begin{proof}[{\it Proof of Proposition \ref{ccoincid}}]
To simplify notations we call
$$\overline{c}_{\phi}= \sup\left\{c(x,y)\,:c\in\mathcal{C}(\phi)\right\},$$
so that we have to prove that $ \overline{c}_{\phi} = c_\phi.$

First we consider the constant sequence $\phi_n:=\phi,\quad\forall n\in\mathbb{N}$. Then $c_\phi\in\mathcal{C}(\phi)$ and we get that $\overline{c}_{\phi}\geq c_\phi$.

Let us prove the converse inequality. 
Let $x,y\in\overline\Omega$, $k>0$ and $\sigma\in H^{x,y}$ such that $L_{\phi}(\sigma)<c_\phi(x,y)+1/k$. 
Let us fix a sequence $\phi_n\rightarrow\phi$ in $L^{q}$ such that $c_{\phi_n}$ converges uniformly to some $c$, we want to prove that $c\leq c_{\phi}$.
From density of simple functions and continuity of $\phi$ we can assume that there exists a finite decomposition $\{t_0, t_1, \cdots t_M\}$ of the interval $[0,1]$ such that $\dot\sigma$ is constant and horizontal on the interval $[t_{i-1},t_i]$; in particular
\begin{equation*}
    L_{\phi_n}(\sigma)=\sum_{i=1}^{M}\int_{t_{i-1}}^{t_i}\phi_n(\sigma(t))|\dot\sigma(t)|_Hdt.
\end{equation*} 
Let us consider a single interval $[t_{i-1}, t_i]$: up to  a change of coordinates, we can also assume that $|\dot\sigma|_H=1$ on this interval. For this reason, in the change of coordinates  $\Psi_i:\mathbb{R}^{2n+1}\rightarrow\mathbb{H}^{n}$, introduced in \eqref{changeofcoordinates}, we can choose 
$\Phi_i(\sigma(t_{i-1}))=(t_{i-1},0)$ so that  
$\Phi_i(\sigma(t_{i}))=(t_{i},0)$, and $$\Phi_i \circ \sigma : [t_{i-1}, t_{i}] \to\mathbb{R}^{2n+1}, \quad (\Phi_i \circ \sigma)(t) = (t, 0),$$
where $\Phi_i:=\Psi_i^{-1}:\mathbb{H}^n\to\mathbb{R}^{2n+1}$.

We now consider, for every $\delta>0$ and for every $i$, cylindrical neighborhoods $C_{i, \delta} = \{(t, \hat e)\in\mathbb{R}^{2n+1}: t \in [t_{i-1}, t_i], |\hat e|_{\mathbb{R}^{2n}} \leq \delta\},$ of the curve $\Phi_i \circ \sigma$, with basis
$S_{i-1} =\{(t_{i-1}, \hat e): |\hat e|_{\mathbb{R}^{2n}} \leq \delta\} .$
For every $\hat e\in \mathbb{R}^{2n},$ with  $|\hat e|_{\mathbb{R}^{2n}}\leq \delta$, we call  $\sigma_e (t) =\Psi_i(t, \hat e)$. By definition 
$$
c_{\phi_n}\Big(\Psi_i(t_{i-1}, \hat e) , \Psi_i(t_i, \hat e)\Big)\leq 
L_{\phi_n}(\sigma_e \circ \theta_i),
$$
where $\theta_i$ is a  change of coordinate which sends $[0,1]$ to $[t_{i-1}, t_{i}]$. 
Note that 
\begin{equation}\label{questa}
L_{\phi_n}(\sigma_e \circ \theta_i) =  L_{\phi_n\circ\Psi_i }( \Phi_i \circ \sigma_e \circ \theta_i) = \int_{t_{i-1}} ^{t_i}  ( \phi_n\circ\Psi_i)(t, \hat e) dt.
\end{equation}
Hence, integrating on $S_{i-1}$  we get
\begin{equation}\label{cpne}
	\int_{S_{i-1}}c_{\phi_n}\Big(\Psi_i(t_{i-1}, \hat e), \Psi_i(t_i, \hat e)\Big) d\mathcal{L}^{2n}(\hat e)\leq \int_{S_{i-1}} \int_{t_{i-1}}^{t_i} ( \phi_n\circ\Psi_i)(t, \hat e) dt d\mathcal{L}^{2n}(\hat e).
\end{equation}
For $ n \to \infty$ using the uniform convergence of $c_{\phi_n}$ to $c$ and the $L^{q}$ convergence of $\phi_n$ to $\phi$ we get
that
\begin{equation*}
    \int_{S_{i-1}} c\Big(\Psi_i(t_{i-1}, \hat e), \Psi_i(t_i, \hat e)\Big)d\mathcal{L}^{2n}(\hat e)\leq\int_{C_i}(\phi\circ\Psi_i) (t, \hat e) d\mathcal{L}^{2n+1} (t, \hat e).
\end{equation*}

Now we divide by the measure of $S_{i-1}$  and pass to the limit as $\delta\rightarrow0^+$.
Using the fact that $c$ is continuous
\begin{equation*}
	\lim_{\delta\to0^+}\frac{1}{d\mathcal{L}^{2n}(S_{i-1})}\int_{S_{i-1}} c\Big(\Psi_i(t_{i-1}, \hat e), \Psi_i(t_i, \hat e)\Big))d\mathcal{L}^{2n}(\hat e)=c\Big(\Psi_i(t_{i-1}, 0), \Psi_i(t_i, 0)\Big) = c(x^{i-1},x^i),
\end{equation*}
where  $x^i = \sigma(t_i)$, 
Analogously  the integral over $C_i =[t_{i-1}, t_i] \times S_{i-1}$ divided by the measure of $S_{i-1} $ converges to the integral on $[t_{i-1}, t_i] $, which is the integral along the curve $\Phi_i\circ\sigma(t)$
\begin{equation*}
	\lim_{\delta\to 0^+}\frac{1}{d\mathcal{L}^{2n}(S_{i-1})}\int_{C_i}(\phi\circ\Psi_i)(t, \hat e) d\mathcal{L}^{2n+1}(t, \hat e)= \int_{t_{i-1}}^{t_i}
  (\phi\circ\Psi_i)(t, 0) dt = 
\int_{t_{i-1}}^{t_i}  \phi(\sigma(t))|\dot\sigma(t)|_Hdt.
\end{equation*}
Then, using \eqref{cpne}, we get that 
\begin{equation*}
 c(x^{i-1},x^{i}) \leq \int_{t_{i-1}}^{t_i}\phi(\sigma(t))|\dot\sigma(t)|_Hdt,\quad \forall i=1,\ldots,M,
\end{equation*}
and then, from \eqref{ineq},
\begin{equation*}
    c(x,y)\leq
 \sum_{i=1}^{M}c(x^{i-1},x^{i}) 
		\leq\sum_{i=1}^{M}\int_{t_{i-1}}^{t_i}\phi(\sigma(t))|\dot\sigma(t)|_Hdt=L_{\phi}(\sigma).
\end{equation*}
This gives 
\begin{equation*}
    c(x,y)\leq c_{\phi}(x,y)+\frac{1}{k}
\end{equation*}
for the choice of $\sigma$ and, since $k$ is arbitrary, it follows that $c(x,y)\leq c_{\phi}(x,y)$.
\end{proof}

Since definition \eqref{barcphi} makes sense also for $L^{q}$ functions, and extends \eqref{ccontinuous}, we will use it as definition of $c_\phi$, for any non-negative function $\phi\in L^{q}(\Omega)$. 

\begin{deff}
If $\phi\in L^{q}(\Omega)$, $\phi\geq0$ then we define
\begin{equation}\label{barcphi}
    c_\phi(x,y)=\sup\left\{c(x,y)\,:c\in\mathcal{C}(\phi)\right\},
\end{equation}
where $\mathcal{C}(\phi)$ has been defined in 
\ref{Setcphi}. 
\end{deff}

Following the same argument as in \cite[Lemma 3.5]{Santambrogio1} we can prove the following result.
\begin{lemma}\label{existapprox}
	Let $q>N$, $\phi\in L^{q}(\Omega)$, $\phi\geq 0$, then there exists a sequence $(\phi_n)_{n\in\mathbb{N}}\subset C(\overline{\Omega},\mathbb{R}_+),\, \phi_n\to\phi\,\mbox{ in }L^{q}(\Omega)$, such that $c_{\phi_n}$ converges to $c_\phi$ in $C(\overline{\Omega}\times\overline{\Omega})$. 
\end{lemma}
 
Finally, as consequence of Remark \ref{30luglio} and Lemma \ref{existapprox}, the function $c_\phi$ is a pseudo distance.

\section{Congested optimal transport problem in $\mathbb{H}^n$}

The scope of this section is to adapt the congested optimal transport problem, originally proposed by Carlier et al. in \cite{Santambrogio1} for the Euclidean setting, to the Heisenberg Group. This adaptation employs the notion of horizontal traffic intensity, which can be viewed as a more abstract version of the horizontal transport density introduced in Section 3. While the latter is defined by integrating over geodesics with extremes in the the space, the former is obtained by integrating over the whole space of horizontal curves.

\subsection{Horizontal traffic plans and traffic intensity}

Let $\Omega\subset\mathbb{H}^{n}$ be a bounded domain with $C^{1,1}$ boundary, describing the area where  the transport problem takes place. Let $\mu,\nu\in\mathcal{P}(\overline\Omega)$ be two probability measures representing the initial and  final states of the system. We introduce the notion of horizontal traffic plan as a  probability measures over the space $H$ of horizontal curves, defined in \eqref{acca}.

A \textit{horizontal traffic plan admissible between} $\mu$ and $\nu$ is a probability measure $Q\in\mathcal{P}(H)$ such that $(e_0)_{\#}Q=\mu$ and $(e_1)_{\#}Q=\nu$, where $e_0$ and $e_1$ are the evaluation maps at times $t=0$ and $t=1$, respectively.

In the Euclidean setting, the convexity of $\Omega$ guarantees the existence of at least one traffic plan. In view of results contained in \cite{Monti2}, in the next remark we have to replace the convexity of $\Omega$ with a different geometric assumption.

\begin{Remark}
Assume that $\overline{\Omega}$ contains the trajectories of geodesics with extremes in the support of $\mu$ and $\nu$ respectively: 
\begin{equation*}
    \mathcal{T}(\mu,\nu):=\left\{S_t(x,y):x\in\textnormal{supp}(\mu), y\in\textnormal{supp}(\nu), t\in[0,1]\right\}\subseteq\overline\Omega,
\end{equation*}
where $S$ is the selection of geodesics fixed in \eqref{31luglio} and $\gamma\in\Pi(\mu,\nu)$, then  \begin{equation}\label{traffpl}
    Q_\gamma:=S_\#\gamma
\end{equation}
is an horizontal  traffic plan between $\mu$ and $\nu$
\end{Remark}

In the sequel we will   denote 
$$\mathcal{Q}_{\mathbb{H}}(\mu,\nu):=\{\text{horizontal traffic plans admissible between $\mu$ and $\nu$}\}.$$

One can associate to  any horizontal traffic plan  $Q\in\mathcal{Q}_{\mathbb{H}}(\mu,\nu)$ a positive Radon measure $i_{Q}$ over $\overline\Omega$, such that
\begin{align}\label{defiQ}
    \int_{\overline{\Omega}} \phi(x) di_{Q}(x):=\int_{H} L_{\phi}(\sigma) d Q(\sigma),\quad \forall \phi \in C(\overline{\Omega},\mathbb{R}_+).
\end{align}
This measure is called the \textit{horizontal traffic intensity}. 
The value $i_{Q}(A)$  of the intensity on a Borel set $A$,  
provides an estimate of how much traffic there is along the horizontal curves in $A$, selected by the traffic assignment $Q$.

\begin{Remark}
The horizontal traffic intensity is a generalization of the horizontal transport density introduced in Section \ref{sectrandens}: indeed, if $\gamma\in\Pi_1(\mu,\nu)$ and $Q_\gamma$ is as in \eqref{traffpl}, then
\begin{equation}\label{transporttraffic}
i_{Q_\gamma}=a_\gamma.
\end{equation}
\end{Remark}

Congestion effects are modeled by a continuous function 
$$g:\mathbb{R}_+\rightarrow\mathbb{R}_+,$$
such that
\begin{enumerate}
    \item $g$ is strictly increasing;
    \item $\lim_{i\to\infty}g(i)=+\infty$.
\end{enumerate}

We refer to such a function as a \textit{congestion function}.

Given $Q\in\mathcal{Q}_{\mathbb{H}}(\mu,\nu)$, we denote by
\begin{equation}\label{congfun}
\phi_{Q}(x):=
\begin{cases}
    g(i_{Q}(x)),\quad &\textnormal{if }i_Q\ll\mathcal{L}^{2n+1},\\
    +\infty,\quad &\textnormal{otherwise},
\end{cases}
\end{equation}
where $i_Q(x)$ is the density of the measure $i_Q$ with respect to the Lebesgue measure $\mathcal{L}^{2n+1}$.  

\begin{Remark}
    The existence of at least one $Q$ such that $i_{Q}\ll\mathcal{L}^{2n+1}$ depends again on $\mu$ and $\nu$. For instance, if either $\mu\ll\mathcal{L}^{2n+1}$ or $\nu\ll\mathcal{L}^{2n+1}$, then the existence of such a $Q$ follow from \eqref{transporttraffic} and Theorem \ref{absolutecont}. 
\end{Remark} 

\subsection{Horizontal Wardrop equilibria}

In this subsection, we introduce the notion of Wardrop equilibrium in the Heisenberg group. This concept was first introduced in the discrete Euclidean setting in \cite{Wardrop} and later formalized in the continuous Euclidean setting in \cite{Santambrogio1}. A Wardrop equilibrium describes the minimization of the cost of transport through two competing minimization problems. First, since the length of curves, weighted by the traffic intensity and the congestion function, represents the cost of transport, an optimal horizontal traffic plan must be concentrated on geodesics with respect to a suitable congested metric. Second, the associated transport plan must solve the Monge-Kantorovich problem associated with this congested metric.

Let $p>1$ and let us suppose that 
\begin{equation}\label{31luglio2}
    \mathcal{Q}_{\mathbb{H}}^p(\mu,\nu):=\left\{Q\in\mathcal{Q}_{\mathbb{H}}(\mu,\nu):i_{Q}\in L^p(\Omega)\right\}\not=\emptyset.
\end{equation}

If for instance $\mu,\nu\in L^p(\Omega)$, then \eqref{31luglio2} holds, see \cite{circelli2024continuous}.

We assume that $\exists a,b\in\mathbb{R}_+$ such that the congestion function satisfies
\begin{equation}\label{growthg}
	ai^{p-1}\leq g(i)\leq b(i^{p-1}+1).
\end{equation}

The $(p-1)$ - growth condition \eqref{growthg} implies that 
\begin{equation}\label{3ottbre}
    \phi_Q:=g\circ i_{Q}\in L^{q}(\Omega),\quad \textnormal{with }q:=\frac{p}{p-1}
\end{equation}
for every $Q\in\mathcal{Q}^p_{\mathbb{H}}(\mu,\nu)$.

The first step is to extend the definition of weighted length of horizontal curves to positive $q$-summable weights.

\begin{teo}\label{prolLxi} Let $Q\in \mathcal{Q}^p_{\mathbb{H}}(\mu,\nu)$. If $\phi\in L^{q}(\Omega),\phi\geq0$, with $q>N$, and $(\phi_n)_{n\in\mathbb{N}}\subset C(\overline{\Omega},\mathbb{R}_+)$ is a sequence such that $\phi_n\rightarrow\phi$ in $L^{q}$, then:
\begin{enumerate}
    \item[(i)] $L_{\phi_n}\to L_\phi$ in $L^1(H,Q)$, where $L_\phi$ is independent of the  approximating sequence $(\phi_n)_{n\in\mathbb{N}}$. 
    \item[(ii)] It holds
    \begin{equation}\label{eglc}
        \int_{\Omega} \phi(x) i_{Q}(x)dx=\int_{H} L_{\phi}(\sigma) dQ(\sigma).
    \end{equation}
    \item[(iii)] It holds
    \begin{equation}\label{ineglc}
		L_{\phi}(\sigma)\geq c_{\phi}(\sigma(0),\sigma(1)),\quad \text{ for }Q-\text{a.e. } \sigma\in H,
    \end{equation}
    where $c_{\phi}$ is defined in \eqref{barcphi}.
\end{enumerate}
\end{teo}

The proof follows from Lemma \ref{existapprox} by passing to the limit, see \cite[Lemma 3.6]{Santambrogio1}.

If we suppose that $p < \frac{N}{N-1}$, then the dual exponent $q > N$. From \eqref{3ottbre}, it follows that both the length $L_{\phi_Q}$, weighted by the traffic intensity and the congestion function, and the associated \textit{congested metric} $c_{\phi_Q}$ are well-defined for every $Q \in \mathcal{Q}_H^p(\mu, \nu)$. See \eqref{barcphi} and Theorem \ref{prolLxi}.


The definition of Wardrop equilibria in $\mathbb{H}^n$ can be given as follows.

\begin{deff}\label{Wardrop}
A horizontal Wardrop equilibrium is a horizontal traffic plan $Q\in\mathcal{Q}^p_{\mathbb{H}}(\mu,\nu)$ such that
\begin{equation}\label{6agosto1}
    Q\big{(}\left\{\sigma\in H: L_{\phi_Q}(\sigma)=c_{\phi_{Q}}(\sigma(0),\sigma(1))\right\}\big{)}=1
\end{equation}
and $\gamma_{Q}:=(e_0,e_1)_{\#}Q\in\Pi(\mu,\nu)$ solves the Monge-Kantorovich problem
\begin{equation}\label{6agosto}
    \inf_{\gamma\in\Pi(\mu,\nu)}\int_{\overline{\Omega}\times\overline{\Omega}}c_{\phi_{Q}}(x,y)d\gamma(x,y).
\end{equation}
\end{deff}

This definition describes an equilibrium in the sense that it is not a priori clear that the two problems can be minimized simultaneously.

We now introduce a convex minimization problem, whose solutions are Wardrop equilibria in $\mathbb{H}^n$. This problem is the Heisenberg analogue of the one proposed in \cite{Santambrogio1}, which was strongly inspired by the discrete case \cite{Beckmann1}.

We call \textit{congested optimal transport problem in} $\mathbb{H}^n$ the following minimization problem
\begin{equation}\label{lepbme1}
	\inf_{Q\in \mathcal{Q}_{\mathbb{H}}^p(\mu,\nu)}\int_{\Omega}G(i_{Q}(x))dx,
\end{equation}
where $G(i):=\int_0^ig(z)dz$.

\begin{teo}
    The minimization problem \eqref{lepbme1} admits solutions.
\end{teo}

\begin{proof}
Due to the convexity of $G$, the result immediately follows from direct methods of calculus of variations see \cite[Theorem 2.10]{Santambrogio1}.
\end{proof}

The following theorem guarantees the existence of horizontal Wardrop equilibria, according to Definition \ref{Wardrop}.

\begin{teo}\label{cwe}
If $1<p<\frac{N}{N-1}$ then $Q\in\mathcal{Q}^p_{\mathbb{H}}(\mu,\nu)$ solves \eqref{lepbme1} if, and only if, it is a horizontal Wardrop equilibrium.
\end{teo}

\begin{proof}[Sketch of the proof]
The proof works as in \cite{Santambrogio1}. Here, we just show why solutions to \eqref{lepbme1} are equilibrium solutions.

As in \cite[Proposition 3.9]{Santambrogio1}, one can also prove that, given a solution $\overline{Q}\in \mathcal{Q}_{\mathbb{H}}^p(\mu,\nu)$ to \eqref{lepbme1}, then
\begin{equation}
	\int_{\Omega} \phi_{\overline{Q}}(x) i_{\overline{Q}}(x)dx=\inf_{\gamma\in\Pi(\mu,\nu)}\int_{\overline{\Omega}\times\overline{\Omega}}c_{\phi_{\overline{Q}}}(x,y)d\gamma(x,y).
\end{equation}

Moreover, if we denote by $\gamma_{\overline Q}:=(e_0,e_1)_{\#}\overline{Q}\in\Pi(\mu,\nu)$, it follows that
\begin{align*}
	\int_{\overline{\Omega}\times \overline{\Omega}} c_{\phi_{\overline{Q}}}(x,y)d\gamma_{\overline{Q}}(x,y)=\int_{H}c_{\phi_{\overline{Q}}}(\sigma(0),\sigma(1))d\overline{Q}(\sigma)\leq\\ \underset{\eqref{ineglc}}{\leq}\int_{H} L_{\phi_{\overline{Q}}}(\sigma) d\overline{Q}(\sigma)= \int_{\Omega} \phi_{\overline{Q}}(x) i_{\overline{Q}}(x)dx=\\= \inf_{\gamma\in \Pi(\mu,\nu)} \int_{\overline{\Omega}\times \overline{\Omega}} c_{\phi_{\overline{Q}}}(x,y)d\gamma(x,y).
\end{align*}
Hence, $\gamma_{\overline Q}$ solves the  Monge-Kantorovich problem associated with the congested metric $c_{\phi_{\overline{Q}}}$
\begin{equation}
	\inf_{\gamma\in \Pi(\mu,\nu)} \int_{\overline{\Omega}\times \overline{\Omega}} c_{\phi_{\overline{Q}}}(x,y)d\gamma(x,y).
\end{equation}
We also observe that
\begin{align*}
	\int_{H} L_{\phi_{\overline{Q}}}(\sigma) d\overline{Q}(\sigma)= \int_{\overline{\Omega}\times\overline{\Omega}} c_{\phi_{\overline{Q}}}(x,y)d\gamma_{\overline Q}(x,y)\\
	=\int_{H}c_{\phi_{\overline{Q}}}(\sigma(0),\sigma(1))d\overline{Q}(\sigma)
\end{align*}
and, since $L_{\phi_{\overline{Q}}}(\sigma)\geq c_{\phi_{\overline{Q}}}(\sigma(0),\sigma(1))$, we get
\[L_{\phi_{\overline{Q}}}(\sigma)=c_{\phi_{\overline{Q}}}(\sigma(0),\sigma(1)) \quad\mbox{ for }  \overline{Q} \mbox{-a.e. } \sigma.\]
\end{proof}

\begin{Remark}
Additionally, if $G$ is strictly convex, then given $Q_{1}$ and $Q_{2}$ which solves \eqref{lepbme1} it follows that $i_{Q_{1}}=i_{Q_{2}}$. In other words, equilibria are not necessarily unique but they all induce the same intensity or, equivalently, the same congested metric.
\end{Remark}


Additional constraints can be imposed on the problem. Specifically, one could minimize the total cost under the constraint $(e_0,e_1)_{\#}Q\in\Pi\subset\Pi(\mu,\nu)$ and, as a particular case, consider $\Pi=\{\overline{\gamma}\}$. 
In this scenario, the set of horizontal traffic plans is defined as:
\begin{equation*}
	\mathcal{Q}_{\mathbb{H}}(\overline{\gamma}):=\left\{Q\in \mathcal{P}(H) \mid (e_0,e_1)_\# Q=\overline{\gamma}\right\}
\end{equation*}
Wardrop equilibria are horizontal traffic plans that belong to:
$$\mathcal{Q}_{\mathbb{H}}^p(\overline{\gamma}):=\left\{Q\in\mathcal{Q}_{\mathbb{H}}(\overline{\gamma}) \mid i_{Q}\in L^p(\Omega)\right\}$$
and satisfy the first condition of Definition \ref{Wardrop}. All previous arguments can be adapted to this case to establish the existence of equilibria.

\section*{Acknowledgments}
M. C. and G.C. are supported by the project PRIN 2022 F4F2LH - CUP J53D23003760006, G.C. is funded by MNESYS PE12 (PE0000006).

The authors thank L. Brasco and S. Rigot for the useful discussions, suggestions and remarks.

\printbibliography

\end{document}